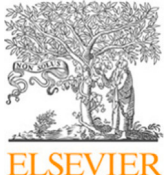
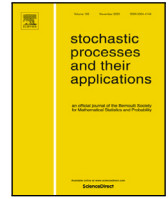

# Large deviations for Markov processes with switching and homogenisation via Hamilton–Jacobi–Bellman equations

Serena Della Corte *, Richard C. Kraaij

*Delft Institute of Applied Mathematics, Delft University of Technology, Mekelweg 4, 2628 CD, Delft, The Netherlands*

A B S T R A C T

We consider the context of molecular motors modelled by a diffusion process driven by the gradient of a weakly periodic potential that depends on an internal degree of freedom. The switch of the internal state, that can freely be interpreted as a molecular switch, is modelled as a Markov jump process that depends on the location of the motor. Rescaling space and time, the limit of the trajectory of the diffusion process homogenises over the periodic potential as well as over the internal degree of freedom. Around the homogenised limit, we prove the large deviation principle of trajectories with a method developed by Feng and Kurtz based on the analysis of an associated Hamilton–Jacobi–Bellman equation with an Hamiltonian that here, as an innovative fact, depends on both position and momenta.

## 1. Introduction

In biochemical and biophysical processes occurring in a cell, an important role is played by several classes of active enzymatic molecules, generally called *motor proteins* or *molecular motors*. These motors are protein molecules that convert chemical energy into mechanical work and motion (see [17–19,33] for more details). In the last decades, such biological phenomena have been largely investigated and this analysis was partly possible due to the contribution of the analysis of particular Markov processes called *"switching Markov processes"* (see for instance [3,4,17,27,34]).

The process that we will consider is such a process. It is a two-component process $(X_t, I_t)$ where the first component $X_t$ is a drift–diffusion process and the second component $I_t$ is a jump process on a finite set. In the context of molecular motors, the spatial component $X_t$ models the location of the motor, for example on a filament, while $I_t$ models the molecular configuration. The two processes together evolve in accordance with the following stochastic differential equation

$$\mathrm{d}X_t = -\nabla \psi(X_t, I_t)\mathrm{d}t + \mathrm{d}B_t,$$
$$\mathbb{P}\Big(I(t+\Delta t) = j \mid I(t) = i, X(t) = x\Big) = r_{ij}(x)\,\Delta t + \mathcal{O}(\Delta t^2) \quad \text{as } \Delta t \to 0,$$

(1.1)

with $\psi \in C^\infty(\mathbb{R}^d \times \{1,\ldots,J\})$, $r_{ij} \in C^\infty(\mathbb{R}^d)$ and $\nabla$ is the gradient with respect to $x$ and $B_t$ is the Brownian motion.

It is clear that the two processes are linked by their rate functions. This means that $I_t$ stays in a first discrete state for a random duration while the diffusion component $X_t$ evolves following a stochastic differential equation with a particular drift. Then, when a switch of the configurational component occurs, the potential $\psi$ changes and therefore $X_t$ diffuses according to a new equation up to another switch of $I_t$ (see Fig. 1 for a typical behaviour of this type of processes). For more details about the construction of such switching hybrid diffusions see [34].

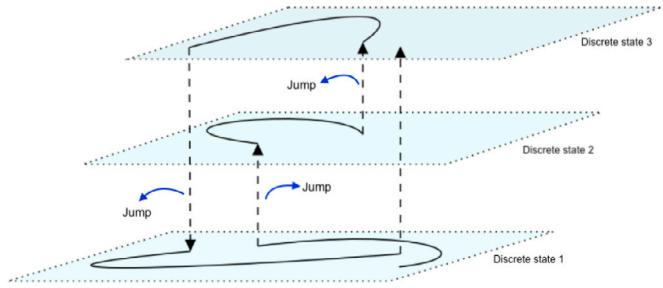

**Fig. 1.** A typical evolution of a process $(X(t), I(t))$.

To allow more flexibility to separate the local dynamics as caused by the internal switching, and macroscopic effects, e.g. modelling the presence of energy molecules in the solution, we will work with $\psi_\varepsilon$ and $r_\varepsilon$ instead of $\psi$ and $r$, and we will typically assume that $\{\psi_\varepsilon, r_\varepsilon\}_\varepsilon$ exhibit a separation of scales. The simplest instance of this separation of scale is that

$$\psi_\varepsilon(x, i) = \psi_1(\varepsilon x, i) + \psi_2(x, i),$$
$$r_\varepsilon(x, i, j) = r_1(\varepsilon x, i, j) + r_2(x, i, j),$$

i.e. $\psi_1$ and $r_1$ model the global macroscopic scale while $\psi_2$ and $r_2$ correspond to the local dynamics. We will also assume that $\psi_2$ is 1–periodic. Moreover, the most general context that we will consider is such that the sequences of functions $\psi_\varepsilon$ and $r_\varepsilon$ are actually given by two functions $\psi \in C^\infty(\mathbb{R}^d \times \mathbb{R}^d \times \{1, \ldots, J\})$ and $r_{ij} \in C^\infty(\mathbb{R}^d \times \mathbb{R}^d)$ as

$$\psi_\varepsilon(x, i) = \psi(\varepsilon x, x, i),$$
$$r_\varepsilon(x, i, j) = r(\varepsilon x, x, i, j). \tag{1.2}$$

The process arising from the stochastic differential equation with $\psi_\varepsilon$ and $r_\varepsilon$ as drift and rate function will be called $(X_t^\varepsilon, I_t^\varepsilon)$. However, we are interested in the macroscopic motion of the molecule. Therefore, we work with the rescaled process or "zoomed out" process that we obtain by scaling in space and time by the positive parameter $\varepsilon > 0$. More precisely, we look at $(Y_t^\varepsilon, \bar{I}_t^\varepsilon) := (\varepsilon X_{\varepsilon^{-1}t}^\varepsilon, I_{\varepsilon^{-1}t}^\varepsilon)$ that evolves according the stochastic differential equation

$$dY_t^\varepsilon = -\nabla \psi \left( Y_t^\varepsilon, \frac{Y_t^\varepsilon}{\varepsilon}, \bar{I}_t \right) dt + \sqrt{\varepsilon}\, dB_t$$
$$\mathbb{P}\left( \bar{I}^\varepsilon(t + \Delta t) = j \mid \bar{I}^\varepsilon(t) = i, Y^\varepsilon(t) = x \right) = \frac{1}{\varepsilon} r_{ij}\left( x, \frac{x}{\varepsilon} \right) \Delta t + \mathcal{O}(\Delta t^2) \quad \text{as } \Delta t \to 0$$

with $\psi$ and $r_{ij}$ given by (1.2), and we are interested in the limit $\varepsilon \to 0$.

Intuitively, looking from far away at the process the periodicity becomes smaller and smaller as $\varepsilon$ decreases and the internal flip rate diverges. Thus, we expect the periodicity and the internal dynamics to homogenise, effectively obtaining a deterministic limit $X_t$. The numerical simulation in Fig. 2 confirms our intuition. It shows sample paths of numerical approximations of a particular switching process $(Y_t^\varepsilon, \bar{I}_t^\varepsilon)$ for various $\varepsilon$. The figure suggests that for small $\varepsilon$, the spatial component $Y_t^\varepsilon$ tends to concentrate around a limiting path that, in the case of the simulated process, is a path with constant velocity.

The aim of this work is to investigate large deviations around such deterministic limit of this kind of process. Showing a large deviations principle, we then will be able to characterise the limit path using a Lagrangian rate function. Indeed, we will show in the main theorem, Theorem 2.8, that there exists a non negative *rate function* $\mathcal{I} : \mathbb{C}_{\mathbb{R}^d}[0, \infty) \to [0, \infty]$ with which $\{Y_t^\varepsilon\}_{\varepsilon > 0}$ satisfies a path-wise large deviation principle in the sense of Definition 2.5 below. Intuitively, it means that

$$\mathbb{P}(Y^\varepsilon \approx x) \sim e^{-\mathcal{I}(x)/\varepsilon} \quad \varepsilon \to 0,$$

with $\mathcal{I}$ written in terms of a Lagrangian function. This means that $Y^\varepsilon$ has a limit path $\tilde{x} \in \mathbb{C}_{\mathbb{R}^d}$ and this limit is the unique minimiser of the rate function $\mathcal{I}$. Moreover, for any path $x \neq \tilde{x}$ such that $\mathcal{I}(x) > 0$, the probability that $Y^\varepsilon$ is close to $x$ is exponentially small in $\varepsilon^{-1}$. In Corollary 2.10, we characterise the minimum of $\mathcal{I}$ by finding a representation of $\partial_t \tilde{x}$ in terms of the drift $\psi$ similar to what one would expect from an averaging principle.

Our work falls into a long tradition of studying the dynamical large deviations around limiting trajectories starting with [15] for small noise diffusions, and [25] for two-scale systems. Then, in the last decades it has been used for the study of different kind of processes (see for instance [13] or [23]). Regarding jump–diffusion, there are very few large deviations results but see for example [24,29]. More recently, in [26], the authors prove large deviations for a class of switching Markov processes and apply their result to examples, including the molecular motors model. Regarding this example, our work diverges from [26] primarily due to two distinctive improvements: the transition from a compact to a non-compact setting and the introduction of the global macroscopic effects in rates with the two components $\psi_1$ and $r_1$. These two facts complicate the proof of large deviations principle. Most important, but without going into details, we need to prove the comparison principle for a spatially inhomogeneous Hamilton–Jacobi–Bellman equation where the two above generalisations introduce non-trivial complications.

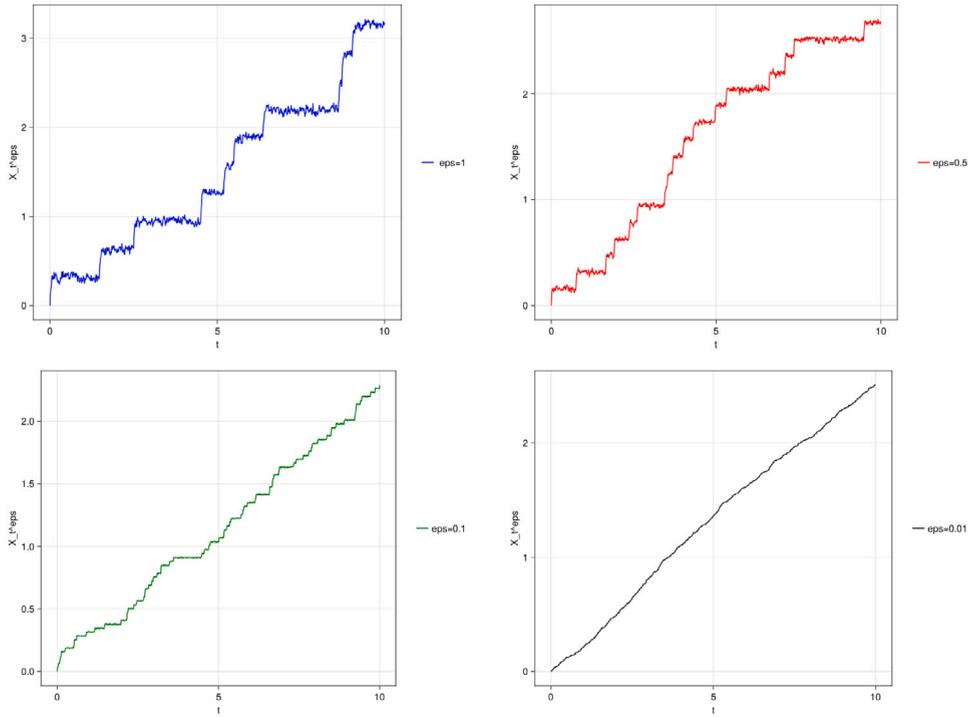

**Fig. 2.** Sample paths of a numerical simulation of the spatial component $Y_t^\varepsilon$ of a switching process for different values of $\varepsilon$. We chose a drift $\psi_1^i$ equal to the periodic part $\psi_2$ for all $i \in \{1, \ldots, 4\}$. We took $\psi_2(\frac{x}{\varepsilon}, i)$ equal to $\sin(x/\varepsilon), \cos(x/\varepsilon), -\sin(x/\varepsilon)$ and $-\cos(x/\varepsilon)$ for $i = 1, 2, 3, 4$ respectively and a rate equal to 1. The jump process switches from a value $i \in \{1, 2, 3\}$ to the value $i + 1$ and from 4 to 1. In this way, the process starts diffusing around a minimum of $\sin(x/\varepsilon)$ ($i = 1$). The first horizontal part of the "stair" corresponds to this evolution. Then, a switch of $\bar{I}_t^\varepsilon$ takes place, so the value of $i$ becomes $i = 2$, and then the spatial component goes to diffuse around a minimum of $\cos(x/\varepsilon)$, that is the second horizontal part, until another switch.

Indeed, we prove the large deviations property using a method due to Feng and Kurtz [14] in which a central role is played by associated Hamilton–Jacobi–Bellmann equations. We will explain in more details the main innovations compared to [26] of our work in Section 6.

The work is organised as follows. We give some preliminary contents and the statement of the large deviations theorem in Section 2. In Section 3 we give an overview of the theory behind the method that we use for showing the large deviations result, proved in Section 4, for the switching process modelling molecular motors. Finally, in Section 5 we are able to extract the main mathematical structures that we use in the previous section and use them in a large deviations result for a more general class of Switching Markov processes.

## 2. General setting and main theorem

In this first section we describe the setting and give some basic notions for the statement of the main theorem. First of all, the following are frequently used notations.

- $C(U, V)$ the space of continuous functions from a set $U \subseteq \mathbb{R}^d$ to a set $V \subseteq \mathbb{R}^d$;
- $C^k(U, V)$, with k integer, the space of k times differentiable functions from $U \subseteq \mathbb{R}^d$ to $V \subseteq \mathbb{R}^d$;
- $C^\infty(U, V)$ the space of infinitely differentiable functions from $U \subseteq \mathbb{R}^d$ to $V \subseteq \mathbb{R}^d$;
- $AC(U, V)$ the space of absolutely continuous functions from $U \subseteq \mathbb{R}^d$ to $V \subseteq \mathbb{R}^d$;
- $\|f\|_\infty = \sup_{x \in E} |f(x)|$;
- $[x]_{\mathbb{Z}^d} = \{y \in \mathbb{R}^d \ : \ x - y \in \mathbb{Z}^d\}$ the equivalence class of $x$ with respect to the relation defined by $\mathbb{Z}^d$;
- $\mathcal{P}(X)$ is the set of probability measures on a space $X$;
- $\mathbb{C}_E[0, \infty)$ is the space of functions defined on $[0, \infty)$ and taking value in a metric space $E$.

### 2.1. Preliminaries

We begin with the definition of the process that we are going to study. It is a two component Markov process $(X_t^\varepsilon, I_t^\varepsilon)$ to which we refer in all the work with *"molecular motors model"* or *"motor proteins model"*.

**Definition 2.1** (*Molecular Motors*). Given an integer $J$, we consider the setting $E = \mathbb{R}^d \times \{1, \ldots, J\}$. For all $i, j$ in $\{1, \ldots, J\}$, let $r_{ij} \in C^\infty(\mathbb{R}^d \times \mathbb{R}^d; [0, \infty))$ denote nonnegative smooth maps, $\psi^i \in C^\infty(\mathbb{R}^d \times \mathbb{R}^d)$ a smooth and $\nabla \psi^i$ its gradient with respect to $x$. We suppose that $\psi^i$ grows at most linearly in the first component and is periodic in the second one. Finally, given the following operator

$$\widetilde{A}_\varepsilon f(x, i) := -\nabla \psi^i(\varepsilon x, x) \cdot \nabla_x f(\cdot, i)(x) + \frac{1}{2} \Delta_x f(\cdot, i)(x) + \sum_{j=1}^J r_{ij}(\varepsilon x, x) [f(x, j) - f(x, i)],$$

we define the $E$-valued Markov process $(X_t^\varepsilon, I_t^\varepsilon)|_{t \geq 0}$ as the solution to the martingale problem corresponding to $\widetilde{A}_\varepsilon$. More precisely, $(X_t^\varepsilon, I_t^\varepsilon)$ is such that for all $f \in D(\widetilde{A}_\varepsilon)$,

$$f(X^\varepsilon(t), I^\varepsilon(t)) - f(X^\varepsilon(0), I^\varepsilon(0)) - \int_0^t \widetilde{A}_\varepsilon f(X^\varepsilon(s), I^\varepsilon(s)) \, ds$$

is a martingale.

**Remark 2.2.** In our case $r_{ij}$ is regular enough that the martingale problem associated to $A_\varepsilon$ is well posed (see [12,31]).

**Remark 2.3.** It is straightforward to see that the above defined process solves the stochastic differential Eq. (1.1) given in the introduction.

We firstly study the above particular model for which we prove the large deviations property. Then, using this model, we lead to a theorem for a general class of processes called *Switching Markov process* (see Section 5).

As mentioned in the introduction, we will work with the rescaled process $(Y_t^\varepsilon, \bar{I}_t^\varepsilon) =$
$= \left( \varepsilon X_{t/\varepsilon}^\varepsilon, I_{t/\varepsilon}^\varepsilon \right)$. Then, by the chain rule, the generator becomes

$$A_\varepsilon f(x, i) = -\nabla \psi^i \left( x, \frac{x}{\varepsilon} \right) \nabla_x f(\cdot, i)(x) + \frac{\varepsilon}{2} \Delta_x f(\cdot, i)(x) + \frac{1}{\varepsilon} \sum_{j=1}^J r_{ij} \left( x, \frac{x}{\varepsilon} \right) [f(x, j) - f(x, i)]. \tag{2.1}$$

We will assume in the main theorem that the matrix $(R_{ij}(x))_{ij} := (\sup_{y \in \mathbb{R}^d} r_{ij}(x, y))_{ij}$ is irreducible. Here we give the rigorous definition.

**Definition 2.4.** We say that a matrix $A = (A_{ij}(x))_{ij \in \{1, \ldots, J\}, x \in \mathbb{R}^d}$ is irreducible if there is no decomposition of $\{1, \ldots, J\}$ into two disjoint sets $\mathcal{J}_1$ and $\mathcal{J}_2$ such that $A_{ij} = 0$ on $\mathbb{R}^d$ whenever $i \in \mathcal{J}_1$ and $j \in \mathcal{J}_2$.

The main goal of this work is to prove that the spatial component $Y^\varepsilon$ of the above Markov process verifies the large deviation principle. Here we give the main definitions which are written down in terms of a general Polish space $\mathcal{X}$ but will later on be applied for e.g. $\mathcal{X} = \mathbb{R}^d$ or $\mathcal{X} = C_{\mathbb{R}^d}([0, \infty))$.

**Definition 2.5.** Let $\{X_\varepsilon\}_{\varepsilon > 0}$ be a sequence of random variables on a Polish space $\mathcal{X}$. Given a function $I : \mathcal{X} \to [0, \infty]$, we say that

  (i) the function $I$ is a *good rate function* if the set $\{x \mid I(x) \leq c\}$ is compact for every $c \geq 0$.
  (ii) the sequence $\{X_\varepsilon\}_{\varepsilon > 0}$ satisfies the *large deviation principle* with good rate function $I$ if for every closed set $A \subseteq \mathcal{X}$, we have

$$\limsup_{\varepsilon \to 0} \varepsilon \log \mathbb{P}[X_\varepsilon \in A] \leq - \inf_{x \in A} I(x)$$

and, for every open set $U \subseteq \mathcal{X}$,

$$\liminf_{\varepsilon \to 0} \varepsilon \log \mathbb{P}[X_\varepsilon \in U] \geq - \inf_{x \in U} I(x).$$

We recall the definition of exponential tightness and the compact containment condition, typical properties that come out in a large deviations context.

**Definition 2.6** (*Exponential Tightness*). A sequence of probability measures $\{P_\varepsilon\}$ on a Polish space $\mathcal{X}$ is said to be *exponentially tight* if for each $a > 0$, there exists a compact set $K_a \subset \mathcal{X}$ such that

$$\limsup_{\varepsilon \to 0} \varepsilon \log P_\varepsilon(K_a^c) \leq -a.$$

In the next we will consider for the above definitions the space $\mathcal{X} = C_{\mathbb{R}^d}[0, \infty)$.

A sequence $\{X_\varepsilon\}$ of $E$-valued random variables is exponentially tight if the corresponding sequence of distributions is exponentially tight.

**Definition 2.7.** We say that the processes $(Z_\varepsilon(t))$ satisfy the *exponential compact containment condition* if for all $T > 0$ and $a > 0$ there is a compact set $K = K(T, a) \subseteq E$ such that

$$\limsup_{\varepsilon \to 0} \varepsilon \log \mathbb{P} \left[ Z_\varepsilon(t) \notin K \text{ for some } t \in [0, T] \right] \leq -a.$$

*2.2. Statement of the main theorem*

Now we state the main theorem in which we prove sufficient conditions for the large deviation property for the spatial component of the switching process defined in Definition 2.1.

**Theorem 2.8** (*Large Deviation for the "Molecular Motors Model"*). *Let $(X_t^\varepsilon, I_t^\varepsilon)$ be the Markov process of Definition 2.1. Suppose that the matrix $R_{ij} = (\sup_{y \in \mathbb{R}^d} r_{ij}(y))_{ij}$ is irreducible. Denote $Y_t^\varepsilon = \varepsilon X_{t/\varepsilon}^\varepsilon$ the rescaled process. Suppose further that at time zero, the family of random variables $\{Y^\varepsilon(0)\}_{\varepsilon > 0}$ satisfies a large deviation principle in $\mathbb{R}^d$ with good rate function $\mathcal{I}_0 : C_{\mathbb{R}^d}[0, \infty) \to [0, \infty]$. Then, the spatial component $\{Y_t^\varepsilon\}$ satisfies a large deviation principle in $C_{\mathbb{R}^d}[0, \infty)$ with good rate function $\mathcal{I} : C_{\mathbb{R}^d}[0, \infty) \to [0, \infty]$ given in the integral form*

$$\mathcal{I}(x) = \begin{cases} \mathcal{I}_0(x(0)) + \int_0^\infty \mathcal{L}(x(t), \dot{x}(t))\, dt & \text{if } x \in AC([0, \infty); \mathbb{R}^d), \\ \infty & \text{else,} \end{cases}$$

*with $\mathcal{L}(x, v) = \sup_p \{p \cdot v - \mathcal{H}(x, p)\}$ the Legendre transform of a Hamiltonian $\mathcal{H}(x, p)$ given in variational form by*

$$\mathcal{H}(x, p) = \sup_{\mu \in \mathcal{P}(E')} \left[ \int_{E'} V_{x,p}(z)\, d\mu(z) - I_{x,p}(\mu) \right], \tag{2.2}$$

*where $E' = \mathbb{T}^d \times \{1, \ldots, J\}$,*

$$V_{x,p}(y, i) = \frac{1}{2} p^2 - p \cdot \nabla_x \psi^i(x, y)$$

*and the map $I_{x,p} : \mathcal{P}(E') \to [0, \infty]$ is the Donsker–Varadhan function, i.e.*

$$I_{x,p}(\mu) = -\inf_\varphi \int_{E'} e^{-\varphi} L_{x,p}(e^\varphi)\, d\mu,$$

*where the infimum is taken over vectors of functions $\varphi(\cdot, i) \in C^2(\mathbb{T}^d)$, and $L_{x,p}$ is the operator defined by*

$$L_{x,p} u(z, i) = \frac{1}{2} \Delta_z u(z, i) + (p - \nabla_x \psi^i(x, z)) \cdot \nabla_z u(z, i) + \sum_{j=1}^J r_{ij}(x, z) [u(z, j) - u(z, i)]. \tag{2.3}$$

**Remark 2.9.** $E'$ captures the periodic behaviour and the internal state. In the homogenisation context described in the introduction, $E'$ is exactly what is being homogenised over while $L_{x,p}$ describes the dynamics on it.

*2.3. Law of large numbers and speed of the limit process*

The following corollary characterises the limit process.

**Corollary 2.10.** *Consider the same assumptions of Theorem 2.8 for the Markov process $(Y_t^\varepsilon, I_t^\varepsilon)$. Then, the spatial component converges almost surely to the path with velocity given by*

$$\partial_t x = \partial_p \mathcal{H}(x, 0) = -\int_{E'} \nabla_x \psi^i(x, y)\, d\mu_x^*(y),$$

*with $\mu_x^*$ the unique stationary measure of the operator $L_{x,0}$ given in (2.3).*

**Proof.** By Theorem A.1 of [26], the spatial component $Y_t^\varepsilon$ converges almost surely to the set of minimisers of the rate function. More precisely,

$$d(Y_t^\varepsilon, \{\mathcal{I} = 0\}) \to 0 \quad \text{a.s. as } \varepsilon \to 0$$

where $\{\mathcal{I} = 0\} = \{x \in C_{\mathbb{R}^d}[0, \infty) : \mathcal{I}(x) = 0\}$. We now prove that this set is actually a singleton and then characterise the unique element.

With this aim, note that by [30, Theorem 23.5], $v$ is a minimiser of $\mathcal{L}$ if and only if $v \in \partial_p \mathcal{H}(x, 0)$. Moreover, by [16, Theorem 4.4.2],

$$\partial_p \mathcal{H}(x, 0) = co \left\{ \bigcup_{\mu \in \mathcal{P}(E')} \partial \left[ \int_{E'} V_{x,0}\, d\mu - I_{x,0}(\mu) \right] \text{ for all } \mu \text{ s.t. } \mathcal{H}(x, 0) = \int_{E'} V_{x,0}\, d\mu - I_{x,0}(\mu) \right\},$$

where with $co$ we refer to the convex hull of a set and $\partial \left[ \int_{E'} V_{x,0}\, d\mu - I_{x,0}(\mu) \right]$ is the differential of the convex functions $\int_{E'} V_{x,p}\, d\mu - I_{x,p}(\mu)$ for $p = 0$.

We know that $\mathcal{H}(x, 0) = 0$ and $V_{x,0}(z) = 0$ for all $z \in E'$. Then, if $\mu_x^*$ is the optimal measure for $\mathcal{H}(x, 0)$, we have that

$$0 = \mathcal{H}(x, 0) = \int_{E'} V_{x,0}(z)\, d\mu_x^* - I_{x,0}(\mu_x^*) = I_{x,0}(\mu_x^*).$$

We can conclude that the optimal $\mu_x^*$ is the unique stationary measure of $L_{x,0}$ (see Proposition Appendix A.1 in the appendix for existence and uniqueness of $\mu_x^*$). We thus find that $\partial_p \mathcal{H}(x,0) = \left\{ \frac{\partial}{\partial p} \mathcal{H}(x,p)|_{p=0} \right\}$ and hence $\mathcal{I}(x) = 0 \iff \partial_t x = \frac{\partial}{\partial p} \mathcal{H}(x,0)$ for almost all $t$ and

$$\partial_t x = \frac{\partial}{\partial p} \mathcal{H}(x,0) = \int_{E'} \frac{\partial V_{x,p}(z)}{\partial p} \bigg|_{p=0} d\mu_x^*(z)$$

$$= -\int_{E'} \nabla_x \psi(x,z) \, d\mu_x^*(z). \quad \square \quad \square$$

**Remark 2.11.** The above corollary confirms the suggestion of Fig. 2 that, when there is no dependence on $x$ in the drift, the spatial component is converging to a path with constant speed. Indeed, for small $\varepsilon$, $Y_t^\varepsilon$ tends to concentrate around a path with a constant velocity $v = \partial_p \mathcal{H}(0)$.

## 3. Connection with Hamilton–Jacobi equations and strategy of proof

We will now present a brief overview of the technical aspects of the Hamilton–Jacobi approach introduced by [14] to the path-space large deviations theory for Markov processes. In the next we give an outline of the main steps of this argument.

Feng and Kurtz in [14, Theorem 5.15] used a variation of the projective limit method [8,9] to prove that the large deviations property can be obtained as a consequence of the large deviations of the finite dimensional distributions and the exponential tightness of the process. By Bryc's theorem and the Markov property, large deviations for the finite dimensional distributions follows by the convergence of the conditional "cumulant function" that forms the semigroup

$$V_\varepsilon(t) f(x) = \varepsilon \log \mathbb{E}\left[ e^{f(X_\varepsilon(t))/\varepsilon} | X(0) = x \right] = \varepsilon \log \int_E e^{f(y)/\varepsilon} \mathbb{P}_\varepsilon(t,x,dy), \quad (3.1)$$

with $\mathbb{P}_\varepsilon(t,x,dy)$ the transition probabilities of $X_t^\varepsilon$. Note that $V_\varepsilon(t) f = \varepsilon \log S_\varepsilon(t) e^{f/\varepsilon}$, where $S_\varepsilon$ is the linear semigroup associated to the generator $A_\varepsilon$. Computing $V_\varepsilon$ and verifying its convergence is usually hard. In analogy to results for linear semigroups and their generators, the convergence of $V_\varepsilon$ follows from the convergence of the *nonlinear generators* $H_\varepsilon$. Formally applying the chain rule to $V_\varepsilon(t)$ in terms of the linear semigroup $S_\varepsilon(t)$ yields the following definition, that can be put on more rigorous grounds as exhibited in [14,22].

**Definition 3.1** (*Nonlinear Generator*). Let $A_\varepsilon$ the generator of a process $X_t^\varepsilon$. The *nonlinear generator* of $X_t^\varepsilon$ is the map defined in the domain $D(H_\varepsilon) = \left\{ f \in C(E) : e^{f(\cdot)/\varepsilon} \in D(A_\varepsilon) \right\}$ by

$$H_\varepsilon f(x) = \varepsilon e^{-f(x)/\varepsilon} A_\varepsilon e^{f(x)/\varepsilon}. \quad (3.2)$$

More precisely, the problem comes down to two steps. First one needs to prove the convergence of the generators $H_\varepsilon \to H$ in a suitable sense. Then, one has to show that the limiting operator generates a semigroup. The sufficient conditions are, using Crandall–Liggett Theorem [6], the *range condition* and the *dissipativity* property.

**Definition 3.2** (*Range Condition*). Let $E$ be an arbitrary metric space and $H : D(H) \subseteq C_b(E) \to C_b(E)$ a nonlinear operator. We say that $H$ satisfies the range condition if

$$\exists \lambda_0 > 0 : D(H) \subset \overline{\mathcal{R}(I - \lambda H)} \quad \text{for all } 0 < \lambda < \lambda_0.$$

**Definition 3.3** (*Dissipative Operator*). We say that an operator $(H, D(H))$ is dissipative if for all $\lambda > 0$,

$$\|(f_1 - \lambda H f_1) - (f_2 - \lambda H f_2)\| \geq \|f_1 - f_2\|$$

for all $f_1, f_2 \in \mathcal{D}(H)$.

The range condition corresponds to existence of classical solutions for the equation $(1 - \lambda H)u = h$. Hence, we can conclude that, in order to prove large deviations, we need the convergence of the nonlinear generators to a dissipative operator $H$ such that the existence of classical solutions for $(1 - \lambda H)u = h$ holds. However, it is well known that the existence of classical solutions is a too strong condition to make this method work in most cases. As observed by Crandall – Lions in [7], the use of *viscosity solutions* allows to overcome this problem. These weak solutions are defined in order to create an extension $\tilde{H}$ of $H$ that automatically satisfies the range condition and that is still dissipative. Below the definitions for both single and multi-valued operators.

**Definition 3.4** (*Sub- and Supersolutions for Single Valued Operators*). Let $H : \mathcal{D}(H) \subseteq C(E) \to C(E)$ be a nonlinear operator. Then for $\lambda > 0$ and $h \in C(E)$, define viscosity sub- and supersolutions of $(1 - \lambda H)u = h$ as follows:

(i) We say that $u : E \to \mathbb{R}$ is a viscosity subsolution if it is bounded and upper semicontinuous and, for every $f \in \mathcal{D}(H)$, there exists a sequence $x_n \in E$ such that

$$\lim_{n \uparrow \infty} u(x_n) - f(x_n) = \sup_x u(x) - f(x),$$

$$\lim_{n \uparrow \infty} u(x_n) - \lambda H f(x_n) - h(x_n) \leq 0.$$

(ii) We say that $v : E \to \mathbb{R}$ is a viscosity supersolution if it is bounded and lower semicontinuous and, for every $f \in \mathcal{D}(H)$, there exists a sequence $x_n \in E$ such that

$$\lim_{n \uparrow \infty} v(x_n) - f(x_n) = \inf_x v(x) - f(x),$$

$$\lim_{n \uparrow \infty} v(x_n) - \lambda H f(x_n) - h(x_n) \geq 0.$$

A function $u \in C(E)$ is called a viscosity solution of $(1 - \lambda H)u = h$ if it is both a viscosity sub- and supersolution.

**Definition 3.5** (*Sub - and Supersolutions for Multivalued Operators*). Let $H \subseteq C(E) \times C(E \times E')$ be a multivalued operator with domain $\mathcal{D}(H) \subseteq C(E)$. Then for $h \in C(E)$ and $\lambda > 0$, define viscosity solutions of $(1 - \lambda H)u = h$ as follows:

(i) $u : E \to \mathbb{R}$ is a viscosity subsolution of $(1 - \lambda H)u = h$ if it is bounded and upper semicontinuous and, for all $(f, g) \in H$, there exists a sequence $(x_n, z_n) \in E \times E'$ such that

$$\lim_{n \uparrow \infty} u(x_n) - f(x_n) = \sup_x u(x) - f(x),$$

$$\limsup_{n \uparrow \infty} u(x_n) - \lambda g(x_n, z_n) - h_1(x_n) \leq 0.$$

(ii) $v : E \to \mathbb{R}$ is a viscosity supersolution of $(1 - \lambda H)u = h$ if it is bounded and lower semicontinuous and, for all $(f, g) \in H$, there exists a sequence $(x_n, z_n) \in E \times E'$ such that

$$\lim_{n \uparrow \infty} v(x_n) - f(x_n) = \inf_x v(x) - f(x),$$

$$\liminf_{n \uparrow \infty} v(x_n) - \lambda g(x_n, z_n) - h_2(x_n) \geq 0.$$

A function $u \in C(E)$ is called a viscosity solution of $(1 - \lambda H)u = h$ if it is both a viscosity sub- and supersolution.

**Remark 3.6.** Consider the definition of subsolutions. Suppose that the test function $f \in \mathcal{D}(H)$ has compact sublevel sets, then instead of working with a sequence $x_n$, there exists $x_0 \in E$ such that

$$u(x_0) - f(x_0) = \sup_x u(x) - f(x),$$

$$u(x_0) - \lambda H f(x_0) - h(x_0) \leq 0.$$

A similar simplification holds in the case of supersolutions.

In the classical context, the range condition, combined with the dissipativity of the operator can be shown to imply unique solvability of the equation $u - \lambda H u = h$. However, for viscosity solutions this argument does not work anymore. The main reason is that viscosity solutions are in general not in the domain of $H$. In order to address this issue, an option can be to suppose that the following *comparison principle* (implying uniqueness) holds.

**Definition 3.7** (*Comparison Principle*). We say that a Hamilton–Jacobi equation $(1 - \lambda H)u = h$ satisfies the comparison principle if for any viscosity subsolution $u$ and viscosity supersolution $v$, $u \leq v$ holds on $E$.

For two operators $H_\dagger, H_\ddagger \subseteq C(E) \times C(E \times E')$, we say that the comparison principle holds if for any viscosity subsolution $u$ of $(1 - \lambda H_\dagger)u = h$ and viscosity supersolution $v$ of $(1 - \lambda H_\ddagger)u = h$, $u \leq v$ holds on $E$.

The theory above was made rigorous in [14,22]. We present the key result in our context and notations.

**Theorem 3.8** (*Adaptation of 7.18 Of [14] to Our Context*). Let $(X^\varepsilon, I^\varepsilon)$ be a Markov process with generator $A_\varepsilon$ and nonlinear generator $H_\varepsilon$ as in Definition 3.1. Consider the semigroup $V_\varepsilon$ defined in (3.1) and suppose the following

(i) $\exists H$  s.t.  $H_\varepsilon$ converges to $H$ in the sense of Definition 4.2,
(ii) $\forall \lambda > 0, h \in C(E)$, the comparison principle holds for $(1 - \lambda H)u = h$,
(iii) $X^\varepsilon$ is exponentially tight,
(iv) $X^\varepsilon(0)$ satisfies large deviation principle with rate function $I_0$.

Then, there exists a unique viscosity solution $R(\lambda)h$ to $(1 - \lambda H)f = h$ and a unique semigroup $V(t) : C_b(E) \to C_b(E)$ such that

1. $\lim_{m \to \infty} R(\frac{t}{m})^m f(x) = V(t)f(x)$ for every $f \in C_b(E), t \geq 0$ and every $x \in E$,
2. $V_\varepsilon$ converges to $V$ in the sense that for any sequence of functions $f_\varepsilon \in C(E)$ and $f \in C(E)$,

$$\text{if} \quad \|f_\varepsilon - f\|_E \xrightarrow{\varepsilon \to 0} 0 \quad \text{then} \quad \|V_\varepsilon(t)f_\varepsilon \to V(t)f\|_E \xrightarrow{\varepsilon \to 0} 0.$$

Moreover, $X_\varepsilon$ satisfies the large deviation principle with rate function $I : \mathbb{C}_E[0, \infty) \to [0, \infty]$ given by

$$I(x) = I_0(x(0)) + \sup_{k \in \mathbb{N}} \sup_{(t_1, \ldots, t_k)} \sum_{i=1}^k I_{t_i - t_{i-1}}(x(t_i) | x(t_{i-1})) \tag{3.3}$$

with $I_t(z|y) = \sup_{f \in C(E)} [f(z) - V(t)f(y)]$ and the supremum above is taken over all finite tuples $t_0 = 0 < t_1 < t_2 < \cdots < t_k$.

## 4. Proof of the main theorem

Using the discussion of the previous section and Theorem 3.8, we can prove Theorem 2.8.

**Proof of Theorem 2.8.** We claim the following five facts:

1. The nonlinear generators $H_\varepsilon f = \varepsilon e^{-f/\varepsilon} A_\varepsilon e^{f/\varepsilon}$ of $Y_t^\varepsilon$ converge to a multivalued operator $H := \{(f, H_{f,\varphi}) : f \in C^2(\mathbb{R}^d),$ $H_{f,\varphi} \in C(\mathbb{R}^d \times E')$ and $\varphi \in C^2(E')\}$,
2. there exists $\tilde{\varphi}$ such that $H_{f,\tilde{\varphi}}(x, z) = H_{\tilde{\varphi}}(x, p, z) = \mathcal{H}(x, p)$ for all $z \in E'$ and $p = \nabla f$,
3. the comparison principle for $(1 - \lambda H)u = h$ holds,
4. $Y^\varepsilon$ verifies the exponential tightness property,
5. the rate function (3.3) can be represented in the following integral form

$$\mathcal{I}(x) = \mathcal{I}_0(x(0)) + \int_0^\infty \mathcal{L}(x(t), \dot{x}(t)) \, dt$$

with $\mathcal{L}(x, v) = \sup_{p \in \mathbb{R}^d} [p \cdot v - \mathcal{H}(x, p)]$ is the Legendre transform of $\mathcal{H}(x, p)$ in (2.2).

We prove the above claims respectively in Propositions 4.3, (4.2), 4.16, 4.22, 4.25 in the following subsections. Then, once the above facts are proved, we can apply Theorem 3.8 and the required large deviation property follows. □

*4.1. The convergence of generators and an eigenvalue problem*

The first step of the proof of large deviations is based on operator convergence. Since the process and its limit do not live in the same space, we cannot work with the usual definition. In the following, we introduce a new definition of limit for functions and multivalued operator on different spaces.

**Definition 4.1.** Let $f_\varepsilon \in C(\mathbb{R}^d \times \{1, \ldots, J\})$ and $f \in C^2(\mathbb{R}^d)$. We say that $LIM f_n = f$ if

$$\|f_\varepsilon - f \circ \eta_\varepsilon\|_{\mathbb{R}^d \times \{1,\ldots,J\}} = \sup_{\mathbb{R}^d \times \{1,\ldots,J\}} |f_\varepsilon - f \circ \eta_\varepsilon| \to 0 \quad \text{as } \varepsilon \to 0,$$

where $\eta_\varepsilon : \mathbb{R}^d \times \{1, \ldots, J\} \to \mathbb{R}^d$ is the projection

$$\eta_\varepsilon(x, i) = x.$$

**Definition 4.2** (*Extended Limit of Multivalued Operators*). Let $H_\varepsilon \subseteq C(\mathbb{R}^d \times \{1, \ldots, J\})$. Define $ex - LIM H_\varepsilon$ as the set

$$ex - LIM H_\varepsilon =$$
$$= \Big\{ (f, H) \in C^2(\mathbb{R}^d) \times C(\mathbb{R}^d \times \mathbb{T}^d \times \{1, \ldots, J\}) | \exists f_\varepsilon \in D(H_\varepsilon) \,:\, f = LIM f_\varepsilon$$
$$\text{and } \|H \circ \eta'_\varepsilon - H_\varepsilon f_\varepsilon\|_{\mathbb{R}^d \times \{1,\ldots,J\}} \to 0 \Big\},$$

where $\eta'_\varepsilon : \mathbb{R}^d \times \{1, \ldots, J\} \to \mathbb{R}^d \times \mathbb{T}^d \times \{1, \ldots, J\}$ is the function $\eta'_\varepsilon(x, i) = \left(x, \left[\frac{x}{\varepsilon}\right]_{\mathbb{Z}^n}, i\right)$.

The following basic example gives the idea of the intuition behind the definitions above.

**Example:** Let $H_\varepsilon f(x, i) = \nabla f(x) + \varepsilon \Delta f(x)$. Then, for every $f \in C^2(\mathbb{R}^d)$ and $\varphi \in C^2(\mathbb{T}^d)$, we define

$$f_\varepsilon(x, i) = f(x) + \varepsilon \varphi\left(\frac{x}{\varepsilon}, i\right) \quad \text{and} \quad H(x, y, i) = \Delta \varphi^i(y).$$

Then, $(f, H) \in ex - LIM H_\varepsilon$.

**Proposition 4.3** (*Convergence of Nonlinear Generator*). Let $E = \mathbb{R}^d \times \{1, \ldots, J\}$ and let $(Y_t^\varepsilon, \bar{I}_t^\varepsilon)$ be the rescaled Markov process with generator $A_\varepsilon$ from (2.1) and let $H_\varepsilon$ be the nonlinear generators defined in Definition 3.1. Then, the multivalued operator $H \subseteq C(\mathbb{R}^d) \times C(\mathbb{R}^d \times \mathbb{T}^d \times \{1, \ldots, J\})$ given by

$$H := \left\{(f, H_{f,\varphi}) \,:\, f \in C^2(\mathbb{R}^d), H_{f,\varphi} \in C(\mathbb{R}^d \times E') \text{ and } \varphi \in C^2(E')\right\},$$

where the images $H_{f,\varphi} : \mathbb{R}^d \times \mathbb{T}^d \times \{1, \ldots, J\} \to \mathbb{R}$ are

$$H_{f,\varphi}(x, y, i) := \frac{1}{2} \Delta_y \varphi^i(y) + \frac{1}{2} \left|\nabla f(x) + \nabla_y \varphi^i(y)\right|^2 - \nabla_x \psi^i(x, y)(\nabla f(x) + \nabla_y \varphi^i(y))$$
$$+ \sum_{j=1}^J r_{ij}(x, y) \left[e^{\varphi(y,j) - \varphi(y,i)} - 1\right],$$

is such that $H \subseteq ex - LIM H_\varepsilon$. Moreover, for all $\varphi$ parametrising the images we have a map $H_\varphi : \mathbb{R}^d \times \mathbb{R}^d \times \mathbb{T}^d \times \{1, \ldots, J\} \to \mathbb{R}$ such that for all $f \in \mathcal{D}(H)$ and any $x \in \mathbb{R}^d$, the images $H_{f,\varphi}$ of $H$ are given by

$$H_{f,\varphi}(x, z') = H_\varphi(x, \nabla f(x), z'), \quad \text{for all } z' \in \mathbb{T}^d \times \{1, \ldots, J\}.$$

**Proof.** We want to prove that $H_\varepsilon$ converges to $H$ in terms of Definition 4.2. With this aim, note that, by the definitions of $A_\varepsilon$ and $H_\varepsilon$, we have

$$H_\varepsilon f(x,i) = \frac{\varepsilon}{2}\Delta_x f^i(x) + \frac{1}{2}|\nabla_x f^i(x)|^2 - \nabla \psi^i\left(x, \frac{x}{\varepsilon}\right)\nabla_x f^i(x)$$
$$+ \sum r_{ij}\left(x, \frac{x}{\varepsilon}\right)\left(e^{(f(x,j)-f(x,i))/\varepsilon} - 1\right).$$

Choosing functions $f_\varepsilon(x,i)$ of the form

$$f_\varepsilon(x,i) = f(x) + \varepsilon\,\varphi\left(\left[\frac{x}{\varepsilon}\right]_{\mathbb{Z}^n}, i\right) = f\circ\eta_\varepsilon(x,i) + \varepsilon\,\varphi\left(\left[\frac{x}{\varepsilon}\right]_{\mathbb{Z}^n}, i\right),$$

we obtain,

$$H_\varepsilon(f_\varepsilon)(x,i) = \frac{\varepsilon}{2}\Delta f(x) + \frac{1}{2}\Delta_y\varphi^i\left(\left[\frac{x}{\varepsilon}\right]_{\mathbb{Z}^n}\right) + \frac{1}{2}\left|\nabla f(x) + \nabla_y\varphi^i\left(\left[\frac{x}{\varepsilon}\right]_{\mathbb{Z}^n}\right)\right|^2$$
$$- \nabla\psi^i(x, x/\varepsilon)\left(\nabla f(x) + \nabla_y\varphi^i\left(\left[\frac{x}{\varepsilon}\right]_{\mathbb{Z}^n}\right)\right) + \sum_{j=1}^J r_{ij}(x, x/\varepsilon)\left[e^{\varphi\left(\left[\frac{x}{\varepsilon}\right]_{\mathbb{Z}^n}, j\right) - \varphi\left(\left[\frac{x}{\varepsilon}\right]_{\mathbb{Z}^n}, i\right)} - 1\right],$$

where $\nabla_y$ and $\Delta_y$ denote the gradient and Laplacian with respect to the variable $y = x/\varepsilon$. We can conclude that

$$\|f\circ\eta_\varepsilon - f_\varepsilon\|_E = \|f(x) - f_\varepsilon(x,i)\|_E = \varepsilon\|\varphi\|_E \to 0 \text{ as } \varepsilon \to 0,$$

and

$$\|H_{f,\varphi}\circ\eta'_\varepsilon - H_\varepsilon f_\varepsilon\|_E = \sup_{(x,i)\in E}\left|H_{f,\varphi}\left(x, \left[\frac{x}{\varepsilon}\right]_{\mathbb{Z}^n}, i\right) - H_\varepsilon f_\varepsilon(x,i)\right|$$
$$= \frac{\varepsilon}{2}\sup_{(x,i)\in E}|\Delta f(x)| \xrightarrow{\varepsilon\to 0} 0. \quad \square\square$$

**Remark 4.4.** Note that for all $f \in D(H)$ the image $H_\varphi$ has the representation

$$H_\varphi(x, p, z) = e^{-\varphi(z)}\left[B_{x,p} + V_{x,p} + R_x\right]e^{\varphi(z)}$$

with $p = \nabla f(x)$ and

$$(B_{x,p}h)(y,i) := \frac{1}{2}\Delta_y h(y,i) + \left(p - \nabla_x\psi^i(x,y)\right)\cdot\nabla_y h(y,i)$$
$$(V_{x,p}h)(y,i) := \left(\frac{1}{2}p^2 - p\cdot\nabla_x\psi^i(x,y)\right)h(y,i),$$
$$(R_x h)(y,i) := \sum_{j=1}^J r_{ij}(x,y)\left[h(y,j) - h(y,i)\right].$$

**Proposition 4.5** (*Existence of an Eigenvalue*). *Let $E' = \mathbb{T}^d \times \{1,\ldots,J\}$ and let $H_\varphi : \mathbb{R}^d \times \mathbb{R}^d \times E' \to R$ the images of $H$ given in Proposition 4.3. Then, for all $p \in \mathbb{R}^d$ there exists an eigenfunction $g_{x,p} \in C^2(E' \times \{1,\ldots,J\})$ with $g^i_{x,p} > 0$ and an eigenvalue $\lambda_{x,p}$ such that*

$$\left[B_{x,p} + V_{x,p} + R_x\right]g_{x,p} = \lambda_{x,p}g_{x,p}.$$

**Proof.** We want to solve the following eigenvalue problem

$$\left[L_{x,p} + R_x\right]g_{x,p} = \lambda_{x,p}g_{x,p} \tag{4.1}$$

where $L_{x,p}$ is a diagonal matrix with $(L_{x,p})_{ii} = (B_{x,p})_i + (V_{x,p})_i$ and $(R_x)_{ij} = r_{ij}$ for $i \neq j$ and $(R_x)_{ii} = \sum_{j=1}^J r_{ij}$.

Guido Sweers showed (see [32]) that there exists $\gamma_{x,p}$ and $g_{x,p} > 0$ such that

$$\left[-L_{x,p} - R_x\right]g_{x,p} = \gamma_{x,p}g_{x,p}$$

when $L_{x,p}$ is a diagonal matrix with $(L_{x,p})_{ii}$ of the type $-\Delta + p\cdot\nabla + c$. Hence, in our case, the equality (4.1) is verified by taking $\lambda_{x,p} = -\gamma_{x,p}$. $\square$

In the next proposition we prove that the images $H_\varphi$ depend only on $x$ and $p$.

**Proposition 4.6.** *Consider the same setting of Proposition 4.5 and let $\mathcal{H}(x,p)$ be the constant depending on $p$ and $x$ given in (2.2). Then, for all $x, p \in \mathbb{R}^d$ there exist a function $\varphi_{x,p} \in C^2(E')$ such that*

$$H_{\varphi_{x,p}}(x, p, z) = \mathcal{H}(x,p) \qquad \text{for all } z \in E'. \tag{4.2}$$

**Proof.** By Proposition 4.5, there exists a function $g_{x,p}$ and a constant $\lambda_{x,p}$ that satisfy the eigenvalue problem for the operator $L_{x,p} + R_{x,p}$ defined in (4.1). By the variational representation established by Donsker and Varadhan in [10], the eigenvalue is equal to the constant $\mathcal{H}(x,p)$ defined in (2.2). Then, equality (4.2) follows from Remark 4.4 and Proposition 4.5 by choosing $\varphi_{x,p} = \log g_{x,p}$. $\square$

### 4.2. Regularity of the Hamiltonian

Before proving the comparison principle, we first show that the map $p \mapsto \mathcal{H}(x,p)$, constructed out of the eigenvalue problem in Propositions 4.5 and 4.6, is convex, coercive and continuous uniformly with respect to $x$.

**Proposition 4.7** (*Convexity and Coercivity of $\mathcal{H}$*). *The map $\mathcal{H} : (x,p) \mapsto \mathcal{H}(x,p)$ in (2.2) is convex in $p$ and coercive in $p$ uniformly with respect to $x$. Precisely,*

$$\lim_{|p| \to \infty} \inf_{x \in K} \mathcal{H}(x,p) = \infty$$

*for every $K$ compact set. Moreover, $\mathcal{H}(x,0) = 0$ for all $x \in \mathbb{R}^d$.*

**Proof.** By Proposition 4.6 the eigenvalue $\mathcal{H}(x,p)$ admits the representation

$$\mathcal{H}(x,p) = -\sup_{g>0} \inf_{z' \in E'} \left\{ \frac{1}{g(z')} \left[ (-B_{x,p} - V_{x,p} - R_x) g \right](z') \right\}$$

$$= \inf_{g>0} \sup_{z' \in E'} \left\{ \frac{1}{g(z')} \left[ (B_{x,p} + V_{x,p} + R_x) g \right](z') \right\}$$

$$= \inf_{\varphi} \sup_{z' \in E'} \left\{ e^{-\varphi(z')} \left[ (B_{x,p} + V_{x,p} + R_x) e^{\varphi} \right](z') \right\} =: \inf_{\varphi} \sup_{z' \in E'} F(x,p,\varphi)(z'),$$

where the map $F$ is given by

$$F(x,p,\varphi)(y,i) = \frac{1}{2} \Delta \varphi^i(y) + \frac{1}{2} |\nabla \varphi^i(y) + p|^2 - \nabla_x \psi^i(x,y)(\nabla \varphi^i(y) + p)$$

$$+ \sum_{j=1}^{J} r_{ij}(x,y) \left[ e^{\varphi^j(y) - \varphi^i(y)} - 1 \right].$$

Note that $F$ is jointly convex in $p$ and $\varphi$. By Proposition 4.6, for every $x, p$ there exists $\varphi_{x,p}$ such that equality holds, i.e. for any $z' \in E'$, we have $\mathcal{H}(x,p) = F(x,p,\varphi_{x,p})(z')$. Therefore, we obtain for $\xi \in [0,1]$ and any $p_1, p_2 \in \mathbb{R}^d$ with corresponding eigenfunctions $e^{\varphi_1}$ and $e^{\varphi_2}$ that

$$\mathcal{H}(x, \xi p_1 + (1-\xi) p_2) = \inf_{\varphi} \sup_{E'} F(x, \xi p_1 + (1-\xi) p_2, \varphi)$$

$$\leq \sup_{E'} F(x, \xi p_1 + (1-\xi) p_2, \xi \varphi_1 + (1-\xi) \varphi_2)$$

$$\leq \sup_{E'} \left[ \xi F(x, p_1, \varphi_1) + (1-\xi) F(x, p_2, \varphi_2) \right]$$

$$\leq \xi \sup_{E'} F(x, p_1, \varphi_1) + (1-\xi) \sup_{E'} F(x, p_2, \varphi_2)$$

$$= \xi \mathcal{H}(x, p_1) + (1-\xi) \mathcal{H}(x, p_2).$$

Regarding coercivity of $\mathcal{H}(x,p)$, we isolate the $p^2$ term in $V_{x,p}$, to obtain

$$\mathcal{H}(x,p) = \frac{p^2}{2} + \inf_{\varphi} \sup_{E'} \left\{ e^{-\varphi} \left[ B_{x,p} - p \cdot \nabla_x \psi + R_x \right] e^{\varphi} \right\}.$$

Any $\varphi \in C^2(E')$ admits a minimum $(y_m, i_m)$ on the compact set $E'$, and with the thereby obtained uniform lower bound

$$\Gamma(x,p,\varphi) = \sup_{E'} \left\{ e^{-\varphi(z_m)} \left[ B_{x,p} - p \cdot \nabla_x \psi + R_x \right] e^{\varphi(z_m)} \right\}$$

$$\geq e^{-\varphi(z_m)} \left[ B_{x,p} - p \cdot \nabla_x \psi + R_x \right] e^{\varphi(z_m)}$$

$$= \underbrace{\frac{1}{2} \Delta_y \varphi(y_m, i_m)}_{\geq 0} + \underbrace{\frac{1}{2} |\nabla_y \varphi(y_m, i_m)|^2}_{= 0} + (p - \nabla_x \psi^{i_m}(x, y_m)) \cdot \underbrace{\nabla_y \varphi(y_m, i_m)}_{= 0}$$

$$+ \sum_{j \neq i} r_{ij}(x, y_m) \underbrace{\left[ e^{\varphi(y_m, j) - \varphi(y_m, i_m)} - 1 \right]}_{\geq 0} - p \cdot \nabla_x \psi^{i_m}(x, y_m) \geq -p \cdot \nabla_x \psi^{i_m}(x, y_m).$$

Using the lower bound $\Gamma(x,p,\varphi) \geq -p \cdot \nabla_x \psi^{i_m}(x, y_m) \geq \inf_{E'} (-p \cdot \nabla_x \psi)$, it follows that, if $K$ is a compact set

$$\inf_{x \in K} \mathcal{H}(x,p) \geq \frac{p^2}{2} - \sup_{x \in K} \sup_{E'} (p \cdot \nabla_x \psi^i(x,y)) \geq \frac{1}{4} p^2 - \sup_{x \in K} \sup_{E'} |\nabla_x \psi^i(x,y)|^2 \xrightarrow{|p| \to \infty} \infty.$$

Regarding $\mathcal{H}(x,0) = 0$, note that $\Gamma(x,0,\varphi) \leq 0$ for all $x$ and $\varphi$. Then we have the first inequality $\mathcal{H}(x,0) \geq \inf_{\varphi} \Gamma(x,0,\varphi) \geq 0$. For the opposite inequality we choose the function $\varphi = (1, \ldots, 1)$ in the representation of $\mathcal{H}$. $\square$

**Proposition 4.8** (*Continuity of $\mathcal{H}$*). *The map $\mathcal{H} : (x,p) \in \mathbb{R}^d \times \mathbb{R}^d \to \mathcal{H}(x,p) \in \mathbb{R}$ is continuous.*

We will prove the continuity of $\mathcal{H}$ by showing that it is lower and upper semicontinuous. For that, we need the following auxiliary results. In particular, for the lower semicontinuity we will make use of the $\Gamma$–convergence in the sense expressed in the following lemma in which we prove that property in a general context. Later, we will use it for $\mathcal{J}(x,p,\theta) = I_{x,p}(\theta)$.

**Lemma 4.9** ($\Gamma$-*Convergence*)**.** *Given two sets $U, V \subseteq \mathbb{R}^d$ and a constant $M \geq 0$ we define $\Theta_{U,V,M}$ as*

$$\Theta_{U,V,M} = \bigcup_{x \in U, p \in V} \{\theta \in \Theta \mid \mathcal{J}(x,p,\theta) \leq M\}.$$

*Let $\mathcal{J} : \mathbb{R}^d \times \mathbb{R}^d \times \Theta \to [0, \infty]$ satisfy the following assumptions:*

(i) *The map $(x, p, \theta) \mapsto \mathcal{J}(x, p, \theta)$ is lower semi-continuous on $\mathbb{R}^d \times \mathbb{R}^d \times \Theta$.*
(ii) *For every $x$ and $p$ fixed and $M \geq 0$, there exist $U_x$ and $U_p$ open and bounded neighbourhoods and a constant $M'$ such that*

$$\mathcal{J}(y, q, \theta) \leq M' \qquad \text{for all } y \in U_x, \ q \in U_p \text{ and } \theta \in \Theta_{\{x\},\{p\},M}.$$

(iii) *For all compact sets $K_1 \subseteq \mathbb{R}^d$ and $K_2 \subseteq \mathbb{R}^d$ and each $M \geq 0$ the collection of functions $\{\mathcal{J}(\cdot,\cdot,\theta)\}_{\theta \in \Theta_{K_1,K_2,M}}$ is equi-continuous.*

*Then if $x_n \to x$ and $p_n \to p$, the functionals $\mathcal{J}_n$ defined by*

$$\mathcal{J}_n(\theta) := \mathcal{J}(x_n, p_n, \theta)$$

*converge in the $\Gamma$-sense to $\mathcal{J}_\infty(\theta) := \mathcal{J}(x,p,\theta)$. That is:*

1. *If $x_n \to x$, $p_n \to p$ and $\theta_n \to \theta$, then $\liminf_{n \to \infty} \mathcal{J}(x_n, p_n, \theta_n) \geq \mathcal{J}(x, p, \theta)$,*
2. *For $x_n \to x$ and $p_n \to p$ and all $\theta \in \Theta$ there are controls $\theta_n \in \Theta$ such that $\theta_n \to \theta$ and $\limsup_{n \to \infty} \mathcal{J}(x_n, p_n, \theta_n) \leq \mathcal{J}(x,p,\theta)$.*

**Proof.** Let $x_n \to x$ and $p_n \to p$ in $\mathbb{R}^d$. If $\theta_n \to \theta$, then by lower semicontinuity (i),

$$\liminf_{n \to \infty} \mathcal{J}(x_n, p_n, \theta_n) \geq \mathcal{J}(x, p, \theta).$$

For the lim-sup bound, let $\theta \in \Theta$. If $\mathcal{J}(x,p,\theta) = \infty$, there is nothing to prove. Thus suppose that $\mathcal{J}(x,p,\theta)$ is finite, i.e., $\theta \in \Theta_{\{x\},\{p\},M}$ for some $M$. Then, by (ii), there exist a bounded neighbourhood $U_x$ of $x$, a bounded neighbourhood $U_p$ of $p$ and a constant $M'$ such that for any $y \in U_x$ and $q \in U_p$,

$$\mathcal{J}(y, q, \theta) \leq M'.$$

Since $x_n \to x$ and $p_n \to p$, the sequences $x_n$ and $p_n$ are, for $n$ large, contained in $U_x$ and $U_p$, respectively. Taking the constant sequence $\theta_n := \theta$, we thus get that $\mathcal{J}(x_n, p_n, \theta_n) \leq M'$ for all $n$ large enough. By (iii), the family of functions $\{\mathcal{J}(\cdot, \cdot, \theta)\}_{\theta \in \Theta_{\bar{U}_x, \bar{U}_p, M'}}$ is equi-continuous, and hence

$$\lim_{n \to \infty} |\mathcal{J}(x_n, p_n, \theta_n) - \mathcal{J}(x, p, \theta)| \leq 0,$$

and the lim-sup bound follows. $\square$

We can now prove that the function $I_{x,p}$ in (2.2) is $\Gamma$-convergent.

**Proposition 4.10** ($\Gamma$-*Convergence of $I_{x,p}$*)**.** *Let $I_{x,p} : \Theta \to [0, \infty]$ the function defined in (2.2). If $x_n \to x$ and $p_n \to p$, the functionals $I_n(\theta) := I_{x_n, p_n}(\theta)$ converge in the $\Gamma$-sense to $I_\infty(\theta) := I_{x,p}(\theta)$.*

**Proof.** Using Lemma 4.9, we only need to prove that $I_{x,p}$ verifies the assumptions.

**Assumption** (i). For any fixed function $u \in \mathcal{D}(L_{x,p})$ such that $u > 0$, the function $(L_{x,p}u/u)$ is continuous. Thus, for any such fixed $u > 0$ it follows that

$$(x, p, \theta) \mapsto \int_{E'} \frac{L_{x,p}u}{u} \, d\theta$$

is continuous on $\mathbb{R}^d \times \mathbb{R}^d \times \Theta$. As a consequence $I(x,p,\theta)$ is lower semicontinuous as the supremum over continuous functions.

**Assumption** (ii). Fix $x$, $p$ and $M \geq 0$. Let $\theta \in \Theta_{x,p,M}$. Then, $I_{x,p}(\theta) = I(x,p,\theta) \leq M$. It follows from [28, Theorem 3] that the density $\frac{d\theta}{dz}$ exists. Moreover, by the same theorem, for all $y$ and $q$ there exist constants $c_1(y,q), c_2(y,q)$ positive, depending continuously on $y$ and $q$, but not on $\theta$, such that

$$I_{y,q}(\theta) \leq c_1(y,q) \int_{E'} |\nabla g_\theta|^2 \, dz + c_2(y,q),$$

where $g_\theta = (d\theta/dz)^{1/2}$ is the square root of the Radon–Nykodym derivative. As the dependence is continuous in $y$ and $q$, we can find two open neighbourhoods, $U \subseteq \mathbb{R}^d$ of $x$ and $V \subseteq \mathbb{R}^d$ of $p$, such that there exist constants $c_1, c_2$ positive, that do not depend on $\theta$, such that for any $y \in U$ and $q \in V$:

$$I_{y,q}(\theta) \leq c_1 \int_{E'} |\nabla g_\theta|^2 \, dz + c_2 := M',$$

obtaining then (ii).

**Assumption** (iii). By the continuity of $r_{ij}$ and $\psi$, assumption (iii) follows from Theorem 4 of [28]. □

The following technical lemma will give us the upper semi-continuity of $\mathcal{H}$.

**Lemma 4.11** (*Lemma 17.30 in [1]*). *Let $\mathcal{X}$ and $\mathcal{Y}$ be two Polish spaces. Let $\phi : \mathcal{X} \to \mathcal{K}(\mathcal{Y})$, where $\mathcal{K}(\mathcal{Y})$ is the space of non-empty compact subsets of $\mathcal{Y}$. Suppose that $\phi$ is upper hemi-continuous, that is if $x_n \to x$ and $y_n \to y$ and $y_n \in \phi(x_n)$, then $y \in \phi(x)$. Let $f : \text{Graph}(\phi) \to \mathbb{R}$ be upper semi-continuous. Then the map $m(x) = \sup_{y \in \phi(x)} f(x,y)$ is upper semi-continuous.*

We can finally prove the continuity of $\mathcal{H}(x,p)$.

**Proof of Proposition 4.8.** We have already showed that $I_{x,p}(\mu)$ is lower semicontinuous and, since $V_p(x,i)$ is continuous and bounded, $\int_{E'} V_{x,p} d\mu$ is continuous. Then, $f(x,p,\mu) := \int_{E'} V_{x,p} d\mu - I_{x,p}(\mu)$ is upper semi-continuous.

Let $x, p \in \mathbb{R}^d$. We know, by Proposition Appendix A.1 in the appendix, that there exists a unique stationary measure $\theta^0_{x,p}$ such that for all $g \in D(L_{x,p})$,

$$\int_{E'} L_{x,p} g(z,i) d\theta^0_{x,p} = 0. \tag{4.3}$$

Let $L^\lambda_{x,p} = \lambda(\lambda - L_{x,p})^{-1} L_{x,p}$ the Hille–Yosida approximation of $L_{x,p}$. Then we have

$$-\int_{E'} \frac{L_{x,p} u}{u} d\theta^0_{x,p} = -\int_{E'} \frac{L^\lambda_{x,p} u}{u} d\theta^0_{x,p} + \int_{E'} \frac{\left(L^\lambda_{x,p} - L_{x,p}\right) u}{u} d\theta^0_{x,p}$$

$$\leq -\int_{E'} \frac{L^\lambda_{x,p} u}{u} d\theta^0_{x,p} + \frac{1}{\inf_{E'} u} \|(L^\lambda_{x,p} - L_{x,p}) u\|_{E'}$$

$$\leq -\int_{E'} L^\lambda_{x,p} \log u\, d\theta^0_{x,p} + o(1).$$

Sending $\lambda \to 0$ and using (4.3) we have that $I_{x,p}(\theta^0_{x,p}) = 0$. Then, $\mathcal{H}(x,p) \geq \int_{E'} V_{x,p} d\theta^0_{x,p}$. Thus, it suffices to restrict the supremum over $\theta \in \phi(x,p)$ where

$$\phi(x,p) := \left\{ \theta \in \mathcal{P}(E') \mid I_{x,p}(\theta) \leq 2\|\Pi(x,p,\cdot)\|_{\mathcal{P}(E')} \right\},$$

where $\|\cdot\|_{\mathcal{P}(E')}$ denotes the supremum norm on $\mathcal{P}(E')$ and we called for simplicity $\Pi(x,p,\theta) = \int_{E'} V_{x,p} d\theta$.

Note that $\|\Pi(x,p,\theta)\|_{\mathcal{P}(E')} < \infty$ by definition of $V_{x,p}$. It follows that

$$\mathcal{H}(x,p) = \sup_{\theta \in \phi(x,p)} \left[ \int_{E'} V_{x,p} d\mu - I_{x,p}(\mu) \right].$$

$\phi(x,p)$ is non-empty as $\theta^0_{x,p} \in \phi(x,p)$ and it is compact because any closed subset of $\mathcal{P}(E')$ is compact. We are left to show that $\phi$ is upper hemi-continuous. Let $(x_n, p_n, \theta_n) \to (x, p, \theta)$ with $\theta_n \in \phi(x_n, p_n)$. We establish that $\theta \in \phi(x,p)$. By the lower semi-continuity of $I$ and the definition of $\phi$ we find

$$I_{x,p}(\theta) \leq \liminf_n I_{x_n,p_n}(\theta_n) \leq \liminf_n 2\|\Pi(x_n,p_n,\cdot)\|_{\mathcal{P}(E')} = 2\|\Pi(x,p,\cdot)\|_{\mathcal{P}(E')}$$

which implies indeed that $\theta \in \phi(x,p)$. Thus, upper semi-continuity follows by an application of Lemma 4.11.

We proceed with proving lower semi-continuity of $\mathcal{H}$. Suppose that $(x_n, p_n) \to (x,p)$, we prove that $\liminf_n \mathcal{H}(x_n, p_n) \geq \mathcal{H}(x,p)$. Let $\theta$ be the measure such that $\mathcal{H}(x,p) = \Pi(x,p,\theta) - I_{x,p}(\theta)$. By Proposition 4.10, there are $\theta_n$ such that $\theta_n \to \theta$ and $\limsup_n I_{x_n,p_n}(\theta_n) \leq I_{x,p}(\theta)$. Moreover, $\Pi(x_n, p_n, \theta_n)$ converges to $\Pi(x,p,\theta)$ by continuity. Therefore,

$$\liminf_{n \to \infty} \mathcal{H}(x_n, p_n) \geq \liminf_{n \to \infty} \left[ \Pi(x_n, p_n, \theta_n) - I_{x_n,p_n}(\theta_n) \right]$$

$$\geq \liminf_{n \to \infty} \Pi(x_n, p_n, \theta_n) - \limsup_{n \to \infty} I_{x_n,p_n}(\theta_n)$$

$$\geq \Pi(x,p,\theta) - I_{x,p}(\theta) = \mathcal{H}(x,p),$$

establishing that $\mathcal{H}$ is lower semi-continuous. □

### 4.3. Comparison principle

In this section we prove the comparison principle for the Hamilton–Jacobi equation in terms of $H$ by relating it to a set of Hamilton–Jacobi equations constructed from $\mathcal{H}$ (Fig. 3). We introduce the operators $H_\dagger, H_\ddagger$ and $H_1, H_2$. In both cases, the new Hamiltonians will serve as natural upper and lower bounds for $\mathbf{H} f(x) = \mathcal{H}(x, \nabla f(x))$ and $H$ respectively, where $\mathcal{H}$ and $H$ are the operators introduced in Propositions 4.6 and 4.3. These new Hamiltonians are defined in terms of a containment function $Y$, which allows us to restrict our analysis to compact sets. Here we give the rigorous definition.

**Definition 4.12** (*Containment Function*). A function $Y : \mathbb{R}^d \to [0, \infty)$ is a containment function for $V_{x,p}$ in (2.2), if $Y \in C^1(\mathbb{R}^d)$ and it is such that

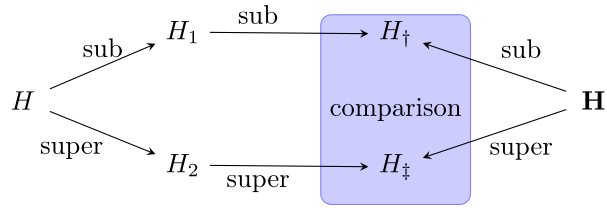

**Fig. 3.** An arrow connecting an operator $A$ with operator $B$ with subscript 'sub' means that viscosity subsolutions of $(1 - \lambda A)f = h$ are also viscosity subsolutions of $(1 - \lambda B)f = h$. Similarly for arrows with a subscript 'super'. The box around the operators $H_†$ and $H_‡$ indicates that the comparison principle holds for subsolutions of $(1 - \lambda H_†)f = h$ and supersolutions of $(1 - \lambda H_‡)f = h$.

- $Y$ has compact sub-level sets, i.e. for every $c \geq 0$ the set $\{x | Y(x) \leq c\}$ is compact ;
- $\sup_{x \in \mathbb{R}^d, z \in E'} V_{x, \nabla Y(x)}(z) < \infty$.

**Lemma 4.13.** *The function $Y(x) = \frac{1}{2} \log \left(1 + |x|^2\right)$ is a containment function for $V_{x,p}$.*

**Proof.** Firstly note that $Y$ has compact sub-level sets. Regarding the second property, by the definition of $V_{x,p}$, we have for every $x \in \mathbb{R}^d$ and $z = (y, i) \in \mathbb{T}^d \times \{1, \ldots, J\}$,

$$V_{x, \nabla Y(x)}(y, i) = \frac{x^2}{2(1 + |x|^2)^2} - \nabla_x \psi^i(x, y) \frac{x}{1 + |x|^2}.$$

Recalling that $\psi$ grows at most linearly in $x$, we can conclude that $\sup_{x,z} V_{x, \nabla Y}(z) < \infty$. □

Using the above lemma we are now able to define the auxiliary operators in terms of $Y$. In the following we will denote by $C_l^\infty(E)$ the set of smooth functions on $E$ that have a lower bound and by $C_u^\infty(E)$ the set of smooth functions on $E$ that have an upper bound.

**Definition 4.14.** Fix $\eta \in (0, 1)$ and given $Y(x) = \frac{1}{2} \log \left(1 + |x|^2\right)$, $C_Y := \sup_{x,z} V_{x, \nabla Y(x)}(z)$ and $\mathbf{H}f(x) = \mathcal{H}(x, \nabla f(x))$, we define

- For $f \in C_l^\infty(E)$,

$$f_†^\eta := (1 - \eta)f + \eta Y,$$
$$H_{†,f}^\eta(x) := (1 - \eta)\mathbf{H}f(x) + \eta C_Y,$$

and set

$$H_† := \left\{ (f_†^\eta, H_{†,f}^\eta) \mid f \in C_l^\infty(E), \eta \in (0, 1) \right\}.$$

- For $f \in C_u^\infty(E)$,

$$f_‡^\eta := (1 + \eta)f - \eta Y,$$
$$H_{‡,f}^\eta(x) := (1 + \eta)\mathbf{H}f(x) - \eta C_Y,$$

and set

$$H_‡ := \left\{ (f_‡^\eta, H_{‡,f}^\eta) \mid f \in C_u^\infty(E), \eta \in (0, 1) \right\}.$$

**Definition 4.15.** Fix $\eta \in (0, 1)$ and given $Y(x) = \frac{1}{2} \log \left(1 + |x|^2\right)$, $C_Y := \sup_{x,z} V_{x, \nabla Y(x)}(z)$ and $\mathbf{H}f(x) = \mathcal{H}(x, \nabla f(x))$, we define

- For $f \in C_l^\infty(E)$, $\varphi \in C^2(E')$, $\eta \in (0, 1)$ set

$$f_1^\eta := (1 - \eta)f + \eta Y,$$
$$H_{1,f,\varphi}^\eta(x, z) := (1 - \eta)H_{f,\varphi}(x, z) + \eta C_Y,$$

and set

$$H_1 := \left\{ (f_1^\eta, H_{1,f,\varphi}^\eta) \mid f \in C_l^\infty(E), \varphi \in C^2(E'), \eta \in (0, 1) \right\}.$$

- For $f \in C_u^\infty(E)$, $\varphi \in C^2(E')$, $\eta \in (0, 1)$ set

$$f_2^\eta := (1 + \eta)f - \eta Y,$$

$$H_{2,f,\varphi}^{\eta}(x,z) := (1+\eta)H_{f,\varphi}(x,z) - \eta C_Y,$$

and set

$$H_2 := \left\{ (f_2^{\eta}, H_{2,f,\varphi}^{\eta}) \mid f \in C_u^{\infty}(E), \varphi \in C^2(E'), \eta \in (0,1) \right\}.$$

We now prove the comparison principle for $f - \lambda H f = h$ based on the results summarised in Fig. 3.

**Theorem 4.16** (*Comparison Principle*). *Let $h \in C_b(E)$ and $\lambda > 0$. Let $u$ and $v$ be, respectively, any subsolution and any supersolution to $(1 - \lambda H)f = h$. Then we have that*

$$\sup_x u(x) - v(x) \leq 0.$$

**Proof.** Fix $h \in C_b(E)$ and $\lambda > 0$. Let $u$ be a viscosity subsolution and $v$ be a viscosity supersolution to $(1 - \lambda H)f = h$. By Fig. 3, the function $u$ is a viscosity subsolution to $(1 - \lambda H_\dagger)f = h$ and $v$ is a viscosity supersolution to $(1 - \lambda H_\ddagger)f = h$. Hence by the comparison principle for $H_\dagger, H_\ddagger$ established in Theorem 4.17 below, $\sup_x u(x) - v(x) \leq 0$, which finishes the proof. □

The rest of this subsection is devoted to establishing Fig. 3. More precisely, we establish Fig. 3 in results 4.17, 4.18, 4.19 and 4.24.

The next theorem contains the comparison principle for $H_\dagger$ and $H_\ddagger$. The proof follows standard ideas that can be found for instance in [2,5]. In order to be able to use both the subsolution and supersolution properties in the estimate of $\sup_x u(x) - v(x)$, we use the following strategy based on the introduction of double variables.

1. First of all, note that the supremum over $x$ of $u(x) - v(x)$ can be replaced, sending $\varepsilon \to 0$, with the supremum over $x$ and $y$ of the double variables function $u(x) - v(y) - (2\varepsilon)^{-1}(x-y)^2$
2. Once the supremum (x,y) is found, we are able to use the sub-super solution properties in the following way:

    - fixing $y$ and optimising over $x$, it can be used in the application of the subsolution property of $u$
    - fixing $x$ and optimising over $y$, it can be used in the application of the supersolution property of $v$.

**Theorem 4.17.** *Let $h \in C_b(E)$ and $\lambda > 0$. Let $u$ be any subsolution to $(1 - \lambda H_\dagger)f = h$ and let $v$ be any supersolution to $(1 - \lambda H_\ddagger)f = h$. Then we have that*

$$\sup_x u(x) - v(x) \leq 0.$$

**Proof.** Following the above steps we define the double variables function

$$\Phi_{\varepsilon,\beta}(x,y) = \frac{u(x)}{1-\beta} - \frac{v(y)}{1+\beta} - \frac{|x-y|^2}{2\varepsilon} - \frac{\beta}{1-\beta}Y(x) - \frac{\beta}{1+\beta}Y(y).$$

Note that the containment function $Y$ is introduced in order to be able to work in a compact set, and the positive constant $\beta$ will allow us to use the convexity of $\mathcal{H}$. Since $\Phi_{\varepsilon,\beta}$ is upper semicontinuous and $\lim_{|x|+|y|\to\infty}\Phi(x,y) = -\infty$, for every $\varepsilon \in (0,1)$ there exists $(x_\varepsilon, y_\varepsilon)$ such that

$$\Phi_{\varepsilon,\beta}(x_\varepsilon, y_\varepsilon) = \sup_{\mathbb{R}^d \times \mathbb{R}^d} \Phi_{\varepsilon,\beta}(x,y). \tag{4.4}$$

Suppose by contradiction that $\delta = u(\tilde{x}) - v(\tilde{x}) > 0$ for some $\tilde{x}$. We choose $\beta$ such that $\frac{2\beta}{(1-\beta)(1+\beta)}Y(\tilde{x}) < \delta/2$ and $\frac{2\beta}{1-\beta^2}(\|h\| + C_Y) < \delta/2$. Then,

$$\Phi_{\varepsilon,\beta}(x_\varepsilon, y_\varepsilon) \geq \Phi_{\varepsilon,\beta}(\tilde{x},\tilde{x}) > \delta - \frac{2\beta}{(1-\beta)(1+\beta)}Y(\tilde{x}) > \frac{\delta}{2} > 0, \tag{4.5}$$

and

$$\frac{\beta}{1-\beta}Y(x_\varepsilon) + \frac{\beta}{1+\beta}Y(y_\varepsilon) \leq \sup\left(\frac{u}{1-\beta}\right) + \sup\left(\frac{-v}{1+\beta}\right) < \infty.$$

Therefore there exists $R_\beta > 0$ such that $x_\varepsilon$ and $y_\varepsilon$ belong to $B(0, R_\beta)$.

Next we observe that by Lemma 3.1 of [5],

$$\frac{|x_\varepsilon - y_\varepsilon|^2}{\varepsilon} \to 0 \qquad \text{as } \varepsilon \to 0^+,$$

and, as a consequence, $|x_\varepsilon - y_\varepsilon| \to 0$ as $\varepsilon \to 0^+$. Define the functions $\varphi_1^{\varepsilon,\beta} \in D(H_\dagger)$ and $\varphi_2^{\varepsilon,\beta} \in D(H_\ddagger)$ by

$$\varphi_1^{\varepsilon,\beta}(x) = (1-\beta)\left[\frac{v(y_\varepsilon)}{1+\beta} + \frac{|x-y_\varepsilon|^2}{2\varepsilon} + \frac{\beta}{1-\beta}Y(x) + \frac{\beta}{1+\beta}Y(y_\varepsilon) + (1-\beta)(x-x_\varepsilon)^2\right]$$

$$\varphi_2^{\varepsilon,\beta}(y) = (1+\beta)\left[\frac{u(x_\varepsilon)}{1-\beta} - \frac{|x_\varepsilon - y|^2}{2\varepsilon} - \frac{\beta}{1-\beta}Y(x_\varepsilon) - \frac{\beta}{1+\beta}Y(y) - (1+\beta)(y-y_\varepsilon)^2\right].$$

Using (4.4), observe that $u - \varphi_1^{\varepsilon,\beta}$ attains its supremum at $x = x_\varepsilon$, and thus

$$\sup_E (u - \varphi_1^{\varepsilon,\beta}) = (u - \varphi_1^{\varepsilon,\beta})(x_\varepsilon).$$

Moreover, by addition of the $(1-\beta)(x-x_\varepsilon)^2$ term, this supremum is the unique optimiser of $u - \varphi_1^{\varepsilon,\beta}$. Then, by the subsolution and supersolution properties, taking into account Remark 3.6,

$$u(x_\varepsilon) - \lambda \left[ (1-\beta)\mathcal{H}\left(x_\varepsilon, \frac{x_\varepsilon - y_\varepsilon}{\varepsilon}\right) + \beta\, C_Y \right] \leq h(x_\varepsilon).$$

With a similar argument for $u_2$ and $\varphi_2^\varepsilon$, we obtain by the supersolution inequality that

$$v(y_\varepsilon) - \lambda \left[ (1+\beta)\mathcal{H}\left(y_\varepsilon, \frac{x_\varepsilon - y_\varepsilon}{\varepsilon}\right) - \beta C_Y \right] \geq h(y_\varepsilon). \tag{4.6}$$

By the coercivity property obtained in Proposition 4.7 on Section 4.2 and by the inequality (4.6), $p_\varepsilon := \frac{x_\varepsilon - y_\varepsilon}{\varepsilon}$ is bounded in $\varepsilon$, allowing us to extract a converging subsequence $p_{\varepsilon_k}$. We conclude that for each $\beta$

$$\liminf_{\varepsilon \to 0} \Phi(x_\varepsilon, y_\varepsilon)$$
$$\leq \liminf_{\varepsilon \to 0} \frac{u(x_\varepsilon)}{1-\beta} - \frac{v(y_\varepsilon)}{1+\beta}$$
$$\leq \liminf_{k \to \infty} \lambda \mathcal{H}(x_{\varepsilon_k}, p_{\varepsilon_k}) + \frac{\beta}{1-\beta} C_Y - \lambda \mathcal{H}(y_{\varepsilon_k}, p_{\varepsilon_k}) + \frac{\beta}{1+\beta} C_Y$$
$$+ \frac{h(x_{\varepsilon_k})}{1-\beta} - \frac{h(y_{\varepsilon_k})}{1+\beta}$$
$$\leq \liminf_{k \to \infty} \lambda \left[ \mathcal{H}(x_{\varepsilon_k}, p_{\varepsilon_k}) - \mathcal{H}(y_{\varepsilon_k}, p_{\varepsilon_k}) \right] + \frac{h(x_{\varepsilon_k}) - h(y_{\varepsilon_k})}{1-\beta^2} + \frac{2\beta}{1-\beta^2} (\|h\| + C_Y)$$
$$\leq \frac{2\beta}{1-\beta^2} (\|h\| + C_Y).$$

As $\beta$ is chosen such that $\frac{2\beta}{1-\beta^2}(\|h\| + C_Y) < \delta/2$, we obtain a contradiction with (4.5), establishing the comparison principle. □

Below, we complete the figure by proving the left-hand side of Fig. 3.

**Lemma 4.18.** *For all $h \in C(\mathbb{R}^d)$ and $\lambda > 0$, viscosity subsolutions of $(1 - \lambda H)f = h$ are viscosity subsolutions of $(1 - \lambda H_1)f = h$, and viscosity supersolutions of $(1 - \lambda H)f = h$ are viscosity supersolutions of $(1 - \lambda H_2)f = h$.*

**Proof.** Fix $\lambda > 0$ and $h \in C_b(E)$. Let $u$ be a subsolution to $(1 - \lambda H)f = h$. We prove it is also a subsolution to $(1 - \lambda H_1)f = h$. Fix $\eta \in (0,1)$, $\varphi \in C^2(E')$ and $f \in C_l^\infty(E)$, so that $(f_1^\eta, H_{1,f,\varphi}^\eta) \in H_1$ with $f_1^\eta$ and $H_{1,f,\varphi}^\eta$ as in Definition 4.15. We will prove that there are $(x_n, z_n)$ such that

$$\lim_n u(x_n) - f_1^\eta(x_n) = \sup_x u(x) - f_1^\eta(x), \tag{4.7}$$

$$\limsup_n u(x_n) - \lambda H_{1,f,\varphi}^\eta(x_n, z_n) - h(x_n) \leq 0. \tag{4.8}$$

Given $M := \eta^{-1} \sup_y u(y) - (1-\eta)f(y) < \infty$, as $u$ is bounded and $f \in C_l^\infty(E)$, we have that the sequence $x_n$ along which the limit in (4.7) is attained, is contained in the compact set $K := \{x | Y(x) \leq M\}$. We define $\gamma : \mathbb{R} \to \mathbb{R}$ as a smooth increasing function such that

$$\gamma(r) = \begin{cases} r & \text{if } r \leq M, \\ M+1 & \text{if } r \geq M+2. \end{cases}$$

Denote by $f_\eta$ the function on $E$ defined by

$$f_\eta(x) = \gamma((1-\eta)f(x) + \eta Y(x)) = \gamma(f_1^\eta(x)).$$

By construction, $f_\eta$ is smooth and constant outside a compact set and thus lies in $\mathcal{D}(H)$. We conclude that $(f_\eta, H_{f_\eta,(1-\eta)\varphi}) \in H$. As $u$ is a viscosity subsolution for $(1 - \lambda H)u = h$, there exist $x_n \in E$ and $z_n \in E'$ with

$$\lim_n u(x_n) - f_\eta(x_n) = \sup_x u(x) - f_\eta(x),$$

$$\limsup_n u(x_n) - \lambda H_{f_\eta,(1-\eta)\varphi}(x_n, z_n) - h(x_n) \leq 0. \tag{4.9}$$

Since $f_1^\eta$ equals $f_\eta$ in $K = \{x|Y(x) \leq M\}$, we also have that

$$\lim_n u(x_n) - f_1^\eta(x_n) = \sup_x u(x) - f_1^\eta(x),$$

establishing (4.7). Convexity of $H_{f,\varphi}(x,z) = H_\varphi(x, \nabla f(x), z)$ in $p$ and $\varphi$ yields for arbitrary $(x,z)$ the elementary estimate

$$\begin{aligned} H_{f_\eta,(1-\eta)\varphi}(x,z) &= H_{(1-\eta)\varphi}(x,(1-\eta)\nabla f(x) + \eta \nabla Y(x), z) \\ &\leq (1-\eta)H_\varphi(x, \nabla f(x), z) + \eta H_0(x, \nabla Y(x), z) \\ &= (1-\eta)H_\varphi(x, \nabla f(x), z) + \eta V_{x,\nabla Y(x)}(z) \\ &\leq H^\eta_{1,f,\varphi}(x,z). \end{aligned}$$

Combining the above inequality with (4.9), we have

$$\limsup_n u(x_n) - \lambda H^\eta_{1,f,\varphi}(x,z) - h(x_n) \leq \limsup_n u(x_n) - \lambda H_{f_\eta,(1-\eta)\varphi}(x_n, z_n) - h(x_n) \leq 0,$$

establishing (4.8). The supersolution statement follows in the same way. □

**Lemma 4.19.** *Fix $\lambda > 0$ and $h \in C_b(E)$.*

(a) *Every subsolution to $(1 - \lambda H_1)f = h$ is also a subsolution to $(1 - \lambda H_\dagger)f = h$.*
(b) *Every supersolution to $(1 - \lambda H_1)f = h$ is also a supersolution to $(1 - \lambda H_\ddagger)f = h$.*

The definition of viscosity solutions, Definition 3.5, is written down in terms of the existence of a sequence of points that maximises $u - f$ or minimises $v - f$. To prove the lemma above, we would like to have the subsolution and supersolution inequalities for any point that maximises or minimises the difference. This is achieved by the following auxiliary lemma.

**Lemma 4.20.** *Fix $\lambda > 0$ and $h \in C_b(E)$.*

(a) *Let $u$ be a subsolution to $(1 - \lambda H_1)f = h$, then for all $(f,g) \in H_1$ and $x_0 \in E$ such that*

$$u(x_0) - f(x_0) = \sup_x u(x) - f(x)$$

*there exists a $z \in E'$ such that*

$$u(x_0) - \lambda g(x_0, z) \leq h(x_0).$$

(b) *Let $v$ be a supersolution to $(1 - \lambda H_2)f = h$, then for all $(f,g) \in H_2$ and $x_0 \in E$ such that*

$$v(x_0) - f(x_0) = \inf_x v(x) - f(x)$$

*there exists a $z \in E'$ such that*

$$v(x_0) - \lambda g(x_0, z) \geq h(x_0).$$

For a proof of the above Lemma see Lemma 5.7 of [23].

**Proof of Lemma 4.19.** We only prove the subsolution statement. Fix $\lambda > 0$ and $h \in C_b(E)$. Let $u$ be a subsolution of $(1 - \lambda H_1)f = h$. We prove that it is also a subsolution of $(1 - \lambda H_\dagger)f = h$. Let $f_1^\eta = (1-\eta)f + \eta Y \in \mathcal{D}(H_1)$ and let $x_0$ be such that

$$u(x_0) - f_1^\eta(x_0) = \sup_x u(x) - f_1^\eta(x).$$

For each $\delta > 0$, since $\mathcal{H}(x,p)$ is a principal eigenvalue for $L_{x,p} + R_x$ (as remarked in Proposition 4.6, there exists a function $g$ such that

$$\mathcal{H}(x_0, \nabla f(x_0)) = g^{-1}\left(L_{x_0, \nabla f(x_0)} + R_{x_0}\right)g. \tag{4.10}$$

As

$$\left(f_1^\eta, (1-\eta)g^{-1}\left(L_{x_0, \nabla f(x_0)} + R_{x_0}\right)g + \eta C_Y\right) \in H_1,$$

we find by the subsolution property of $u$ and that there exists $z$ such that

$$\begin{aligned} h(x_0) &\geq u(x_0) - \lambda\left((1-\eta)g^{-1}\left(L_{x_0, \nabla f(x_0)} + R_{x_0}\right)g + \eta C_Y\right) \\ &= u(x_0) - \lambda\left((1-\eta)\mathcal{H}(x_0, \nabla f(x_0)) + \eta C_Y\right) \end{aligned}$$

where the second inequality follows by (4.10) and it establishes that $u$ is a subsolution for $(1 - \lambda H_\dagger)f = h$. □

We conclude this subsection proving the right part of Fig. 3.

**Proposition 4.21.** *Let the map $\mathcal{H} : \mathbb{R}^d \times \mathbb{R}^d \to \mathbb{R}$ be the eigenvalue (2.2) and let $\mathbf{H} : \mathcal{D}(\mathbf{H}) \subseteq C^1(\mathbb{R}^d) \to C(\mathbb{R}^d)$ be the operator $\mathbf{H}f(x) := \mathcal{H}(x, \nabla f(x))$. Then, for all $\lambda > 0$ and $h \in C(\mathbb{R}^d)$, every viscosity subsolution of $(1 - \lambda \mathbf{H})f = h$ is also a viscosity subsolutions of $(1 - \lambda \mathbf{H}_\dagger)f = h$ and every viscosity supersolution of $(1 - \lambda \mathbf{H})f = h$ is also a viscosity supersolution of $(1 - \lambda \mathbf{H}_\ddagger)f = h$.*

**Proof.** Fix $\lambda > 0$ and $h \in C_b(E)$. Let $u$ be a subsolution to $(1 - \lambda \mathbf{H})f = h$. We prove it is also a subsolution to $(1 - \lambda H_\dagger)f = h$. Fix $\eta > 0$ and $f \in C_\ell^\infty(E)$ and let $(f_\dagger^\varepsilon, H_{\dagger,f}^\eta) \in H_\dagger$ as in Definition 4.14. We will prove that

$$\left(u - f_\dagger^\eta\right)(x) = \sup_{x \in E}\left(u(x) - f_\dagger^\eta(x)\right),$$

implies

$$u(x) - \lambda H_{\dagger,f}^\eta(x) - h(x) \leq 0. \tag{4.11}$$

As $u$ is a viscosity subsolution for $(1 - \lambda \mathbf{H})f = h$ and $f_\dagger^\eta \in D(\mathbf{H})$, if

$$\left(u - f_\dagger^\eta\right)(x) = \sup_x \left(u(x) - f_\dagger^\eta(x)\right),$$

then,

$$u(x) - \lambda \mathbf{H} f_\dagger^\eta(x) - h(x) \leq 0. \tag{4.12}$$

Convexity of $p \mapsto \mathcal{H}(x,p)$ yields the estimate

$$\mathbf{H} f_\eta(x) = \mathcal{H}(x, \nabla f_\eta(x))$$
$$\leq (1 - \eta)\mathcal{H}(x, \nabla f(x)) + \eta \mathcal{H}(x, \nabla Y(x))$$
$$\leq (1 - \eta)\mathcal{H}(x, \nabla f(x)) + \eta C_Y = H_{\dagger,f}^\eta(x).$$

Combining this inequality with (4.12), we have

$$u(x) - \lambda H_{\dagger,f}^\eta(x) - h(x) \leq u(x) - \lambda \mathbf{H} f_\dagger^\eta(x) - h(x) \leq 0,$$

establishing (4.11). The supersolution statement follows in a similar way. □

### 4.4. Exponential tightness

To establish exponential tightness, we first note that by [14, Corollary 4.19] it suffices to establish the exponential compact containment condition. This is the content of the next proposition.

**Proposition 4.22.** *For all $K \subset E$ compact, $T > 0$ and $a > 0$ there is a compact set $\hat{K}_{K,T,a} \subset E$ such that*

$$\limsup_{\varepsilon \to 0} \varepsilon \log \mathbb{P}\left[\bigcup_{t \in [0,T]} \{Y^\varepsilon(t) \notin \hat{K}_{K,T,a}\} \neq \emptyset\right] \leq \max\{-a, \limsup_{\varepsilon \to 0} \varepsilon \log \mathbb{P}(X_\varepsilon(0) \notin K)\}. \tag{4.13}$$

**Remark 4.23.** Note that, since $Y^\varepsilon(0)$ satisfies the large deviations principle by assumption, inequality (4.13) gives the searched compact containment condition.

**Proof of Proposition 4.22.** First of all let us consider $\varphi \equiv 0$. Note that, by Lemma 4.13, we have $\sup_{x,z} H_0(x, \nabla Y, z) = \sup_{x,z} V_{x,\nabla Y(x)}(z) \leq C_Y$. Choose $\beta > 0$ such that $TC_Y - \beta \leq -a$. Since $Y$ is continuous, there is some $c$ such that the set $G := \{x \mid Y(x) < c + \beta\}$ is non empty. Note that $G$ is open and let $\overline{G}$ be the closure of $G$. Then, $\overline{G}$ is compact. Let $f(x) := \iota \circ Y$ where $\iota$ is some smooth increasing function such that

$$\iota(r) = \begin{cases} r & \text{if } r \leq \beta + c, \\ 2\beta + c & \text{if } r \geq \beta + c + 2. \end{cases}$$

It follows that $\iota \circ Y$ equals $Y$ on $\overline{G}$ and is constant outside of a compact set. Set $f_\varepsilon = f \circ \eta_\varepsilon$, $g_\varepsilon = H_\varepsilon f_\varepsilon$ and $g = H_{f,\varphi}$. Note that $g(x,z) = H_\varphi(x, \nabla Y(x), z)$ if $x \in \overline{G}$. Therefore, we have $\sup_{x \in \overline{G}, z \in E'} g(x,z) \leq C_Y$. Let $\tau$ be the stopping time $\tau := \inf\left\{t \geq 0 \mid Y^\varepsilon(t) \notin \overline{G}\right\}$ and let

$$M_\varepsilon(t) := \exp\left\{\frac{1}{\varepsilon}\left(f(Y^\varepsilon(t)) - f(Y^\varepsilon(0)) - \int_0^t g_\varepsilon(Y^\varepsilon(s), \overline{I}_\varepsilon(t))\mathrm{d}s\right)\right\}.$$

By construction $M_\varepsilon$ is a martingale. Let $K \subset E$ be compact. We have

$$\mathbb{P}\left[\bigcup_{t \in [0,T]} \left\{Y^\varepsilon(t) \notin \overline{G}\right\} \neq \emptyset\right]$$

$$\leq \mathbb{P}\left(Y^\varepsilon(0) \in K, \bigcup_{t \in [0,T]} \left\{Y^\varepsilon(t) \notin \overline{G}\right\}\right) + \mathbb{P}\left(Y^\varepsilon(0) \notin K\right)$$

$$= \mathbb{E}\left[\mathbb{1}_{\{Y^\varepsilon(0) \in K\}} \mathbb{1}_{\left\{\bigcup_{t \in [0,T]}\left\{Y^\varepsilon(t) \notin \overline{G}\right\}\right\}} M_\varepsilon(\tau) M_\varepsilon(\tau)^{-1}\right] + \mathbb{P}\left(Y^\varepsilon(0) \notin K\right)$$

$$\leq \exp\left\{-\frac{1}{\varepsilon}\left(\inf_{y_1 \notin \overline{G}} f(y_1) - f(Y^\varepsilon(0))\right.\right.$$
$$\left.\left. -T \sup_{y_2 \in \overline{G}, i \in \{1,\ldots,J\}} g_\varepsilon(y_2, i)\right)\right\}$$
$$\times \mathbb{E}\left[\mathbb{1}_{\{Y^\varepsilon(0) \in K\}} \mathbb{1}_{\left\{\bigcup_{t \in [0,T]} \{Y^\varepsilon(t) \notin \overline{G}\}\right\}} M_\varepsilon(\tau)\right] + \mathbb{P}(Y^\varepsilon(0) \notin K).$$

Since $\sup_{x \in \overline{G}, z \in E'} g(x, z) \leq C_{Y,\varphi}$ and $g$ is the limit of $g_\varepsilon$ for $\varepsilon \to 0$ in the sense of Definition 4.2, we obtain that the term in the exponential is bounded by $\frac{1}{\varepsilon}\left(C_Y T - \beta\right) \leq -\frac{1}{\varepsilon} a$ for sufficiently small $\varepsilon$. The expectation is bounded by 1 due to the martingale property of $M_\varepsilon(\tau)$. We can conclude that

$$\limsup_{\varepsilon \to 0} \varepsilon \log \mathbb{P}\left[\bigcup_{t \in [0,T]} \{Y^\varepsilon(t) \notin K_{T,a}\} \neq \emptyset\right] \leq \max\{-a, \limsup_{\varepsilon \to 0} \varepsilon \log \mathbb{P}(Y^\varepsilon(0) \notin K)\}$$

where $\hat{K}_{K,T,a} = \overline{G}$. □

*4.5. Action-integral representation of the rate function*

In this section we establish a representation of the rate function as an integral of a Lagrangian function $\mathcal{L}$. We refer to this representation as the "action-integral representation" of the rate function $\mathcal{I}$. We argue on basis of Section 8 of [14] for which we need to check the following two conditions.

**Lemma 4.24.** *Let* $\mathcal{H} : \mathbb{R}^d \times \mathbb{R}^d \to \mathbb{R}$ *be the map given in* (2.2) *and* $\mathbf{H} : \mathcal{D}(\mathbf{H}) \subseteq C^1(\mathbb{R}^d) \to C(\mathbb{R}^d)$ *the operator* $\mathbf{H}f(x) := \mathcal{H}(x, \nabla f(x))$. *Then:*

(i) *The Legendre–Fenchel transform* $\mathcal{L}(x, v) := \sup_{p \in \mathbb{R}^d}(p \cdot v - \mathcal{H}(x, p))$ *and the operator* $\mathbf{H}$ *satisfy Conditions 8.9, 8.10 and 8.11 of* [14].
(ii) *For all* $\lambda > 0$ *and* $h \in C(\mathbb{R}^d)$, *the comparison principle holds for* $(1 - \lambda \mathbf{H})u = h$.

**Proof.** To prove the first aim, we will show that following items (a), (b) and (c) imply Condition 8.9, 8.10 and 8.11 of [14]. Then, the proof of (a), (b), (c) is shown in [26, Proposition 6.1].

(a) The function $\mathcal{L} : \mathbb{R}^d \times \mathbb{R}^d \to [0, \infty]$ is lower semicontinuous and for every $C \geq 0$, the level set $\{(x, v) \in \mathbb{R}^d \times \mathbb{R}^d : \mathcal{L}(x, v) \leq C\}$ is relatively compact in $\mathbb{R}^d \times \mathbb{R}^d$.
(b) For all $f \in \mathcal{D}(H)$ there exists a right continuous, nondecreasing function $\psi_f : [0, \infty) \to [0, \infty)$ such that for all $(x, v) \in \mathbb{R}^d \times \mathbb{R}^d$,

$$|\nabla f(x) \cdot v| \leq \psi_f(\mathcal{L}(x, v)) \qquad \text{and} \qquad \lim_{r \to \infty} \frac{\psi_f(r)}{r} = 0.$$

(c) For each $x_0 \in \mathbb{R}^d$ and every $f \in \mathcal{D}(\mathbf{H})$, there exists an absolutely continuous path $x : [0, \infty) \to \mathbb{R}^d$ such that $x_0 = x(0)$ and

$$\int_0^t \mathcal{H}(x(s), \nabla f(x(s))) \, ds = \int_0^t [\nabla f(x(s)) \cdot \dot{x}(s) - \mathcal{L}(x(s), \dot{x}(s))] \, ds.$$

Then regarding Condition 8.9, the operator $Af(x, v) := \nabla f(x) \cdot v$ on the domain $\mathcal{D}(A) = \mathcal{D}(H)$ satisfies (1). For (2), we can choose $\Gamma = \mathbb{T}^d \times \mathbb{R}^d$, and for $x_0 \in \mathbb{T}^d$, take the pair $(x, \lambda)$ with $x(t) = x_0$ and $\lambda(dv \times dt) = \delta_0(dv) \times dt$. Part (3) is a consequence of (a) from above. Part (4) can be verified as follows. Let $Y$ the containment function used in Definition 4.14 and note that the sub-level sets of $Y$ are compact. Let $\gamma \in \mathcal{AC}$ with $\gamma(0) \in K$ and such that the control

$$\int_0^T \mathcal{L}(\gamma(s), \dot{\gamma}(s)) \, ds \leq M$$

implies $\gamma(t) \in \hat{K}$ for all $t \leq T$, with $\hat{K}$ compact. Then,

$$Y(\gamma(t)) = Y(\gamma(0)) + \int_0^t \langle \nabla Y(\gamma(s)), \dot{\gamma}(s) \rangle \, ds$$
$$\leq Y(\gamma(0)) + \int_0^t \mathcal{L}(\gamma(s), \dot{\gamma}(s)) + \mathcal{H}(\gamma(s), \nabla Y(\gamma(s))) \, ds$$
$$\leq \sup_{y \in K} Y(y) + M + \int_0^T \sup_z \mathcal{H}(z, \nabla Y(z)) \, ds$$
$$:= C < \infty.$$

Hence, we can take $\hat{K} = \{z \in \mathbb{R}^d | Y(z) \leq C\}$.

Part (5) is implied by (b) from above. Condition 8.10 is implied by Condition 8.11 and the fact that $\mathbf{H}1 = 0$, by Proposition 4.7 (see Remark 8.12 (e) in [14]). Finally, Condition 8.11 is implied by (c) above, with the control $\lambda(dv \times dt) = \delta_{\dot{x}(t)}(dv) \times dt$.

The comparison principle for $\mathbf{H}$ follows from Proposition 4.21 and Theorem 4.17. □

In the following, we prove the integral representation of the rate function. Firstly, let us recall that Theorem 3.8 gives the existence of a semigroup $V(t)$ and a family of functions $R(\lambda)$ and let $\mathbf{V}(t) : C(\mathbb{R}^d) \to C(\mathbb{R}^d)$ be the Nisio semigroup with cost function $\mathcal{L}$, that is

$$\mathbf{V}(t)f(x) = \sup_{\substack{\gamma \in AC_{\mathbb{R}^d}[0,\infty) \\ \gamma(0)=x}} \left[ f(\gamma(t)) - \int_0^t \mathcal{L}(\gamma(t), \dot{\gamma}(s))\,ds \right].$$

Let $\mathbf{R}(\lambda)h$ be the operator given by

$$\mathbf{R}(\lambda)h(x) = \sup_{\substack{\gamma \in AC \\ \gamma(0)=x}} \int_0^\infty \lambda^{-1} e^{-\lambda^{-1}t} \left[ h(\gamma(t)) - \int_0^t \mathcal{L}(\gamma(s), \dot{\gamma}(s))\,ds \right] dt.$$

The proof of the result below is based on the following four main steps.

- Fig. 3 on Section 4.3 shows that $R(\lambda)h$ is the unique function that is a sub- and supersolution to the equations $(1 - \lambda H_\dagger)f = h$ and $(1 - \lambda H_\ddagger)f = h$ respectively.
- $\mathbf{R}(\lambda)h$ has been proven to be the unique viscosity solution to $(1 - \lambda \mathbf{H})f = h$. Then, again by Fig. 3, we must have $R(\lambda)h = \mathbf{R}(\lambda)h$.
- Starting from the equality of resolvents we work to an equality for the semigroups $V(t)$ and $\mathbf{V}(t)$.
- Recalling that the rate function in Theorem 2.8 is given by,

$$I(x) = I_0(x(0)) + \sup_{k \in \mathbb{N}} \sup_{(t_1,\ldots,t_k)} \sum_{i=1}^k I_{t_i - t_{i-1}}(x(t_i)|x(t_{i-1}))$$

with $I_t(z|y) = \sup_{f \in C(E)}[f(z) - V(t)f(y)]$, it is not difficult to realise that, if $V(t) = \mathbf{V}(t)$, it follows that $I_t(y|z) = \inf_{\substack{\gamma : \gamma(0)=z, \\ \gamma(t)=y}} \int_0^t \mathcal{L}(\gamma(s), \dot{\gamma}(s))\,ds$.

**Theorem 4.25** (*Integral Representation of the Rate Function*). *The rate function of* Theorem 2.8 *has the following representation*

$$\mathcal{I}(x) = \begin{cases} \mathcal{I}_0(x(0)) + \int_0^\infty \mathcal{L}(x(t), \dot{x}(t))\,dt & \text{if } x \in AC([0,\infty); \mathbb{R}^d), \\ \infty & \text{else}, \end{cases}$$

*where* $\mathcal{L}(x, v) = \sup_{p \in \mathbb{R}^d}[p \cdot v - \mathcal{H}(x, p)]$ *is the Legendre transform of* $\mathcal{H}$.

**Proof.** Following the above mentioned steps, we recall that, as stated by Theorem 3.8, there exists a family of operators $R(\lambda) : C_b(\mathbb{R}^d) \to C_b(\mathbb{R}^d)$, such that for $\lambda > 0$ and $h \in C_b(\mathbb{R}^d)$, the function $R(\lambda)h$ is the unique function that is a viscosity solution to $(1 - \lambda H)f = h$ and such that

$$\lim_{m \to \infty} \left\| R\left(\frac{t}{m}\right)^m f - V(t)f \right\| = 0 \qquad \text{for all } f \text{ in a dense set } D \subseteq C_b(\mathbb{R}^d). \tag{4.14}$$

See also [21, Theorem 7.10] or [14, Theorem 7.17] for the construction of the operators $R(\lambda)$. By [23, Proposition 6.1] (or [14, Chapter 8]), $\mathbf{R}(\lambda)$ is the unique viscosity solution to $(1 - \lambda \mathbf{H})f = h$. Then, Fig. 3 on Section 4.3 shows that it must equal $R(\lambda)h$. Moreover, we find by [14, Lemma 8.18] (whose assumptions are implied by Lemma 4.24 above) that for all $f \in C_b(\mathbb{R}^d)$ and $x \in \mathbb{R}^d$

$$\lim_{m \to \infty} \mathbf{R}\left(\frac{t}{m}\right)^m f(x) = \mathbf{V}(t)f(x). \tag{4.15}$$

We conclude from (4.14) and (4.15) that $V(t)f = \mathbf{V}(t)f$ for all $t$ and $f \in D$. Now recall that $D$ is sequentially strictly dense so that equality for all $f \in C_b(\mathbb{R}^d)$ follows if $V(t)$ and $\mathbf{V}(t)$ are sequentially continuous. The first statement follows by Theorems [20, Theorem 7.10] and [21, Theorem 6.1]. The second statement follows by [14, Lemma 8.22]. We conclude that $V(t)f = \mathbf{V}(t)f$ for all $f \in C_b(\mathbb{R}^d)$ and $t \geq 0$. Using Theorem 8.14 of [14] and the convexity of $v \mapsto \mathcal{L}(x, v)$ we get the integral representation. □

## 5. A more general theorem

Analysing the proofs in the previous sections, we can state the following facts:

- In the proof of large deviations principle, the main steps are:
  1. Convergence of the nonlinear operators $H_\varepsilon$ to a multivalued operator $H$,
  2. comparison principle for $(1 - \lambda H)f = h$.
- The existence of an eigenvalue $\mathcal{H}(x, p)$ and its convexity, coercivity and continuity are crucial for our approach to comparison principle and
  - the arguments for existence, convexity and coercivity (proofs of Propositions 4.1 and 4.7) are based on the fact that $\mathcal{H}(x, p)$ is the eigenvalue of an operator of the type $B_{x,p} + V_{x,p} + R_x$ with the three operators that verify particular properties such as coercivity and the maximum principle,

- to show the continuity of $\mathcal{H}$ the representation (2.2) is needed. In particular, some properties of $V$ and $I$, like $\Gamma$-convergence, are necessary.

The above observations allow for a straightforward generalisation in Theorem 5.3 and justify the assumptions of the next subsection. In this section we indeed prove the large deviation principle for a general switching Markov process. In particular, we will study the Markov process $(Y_t^\varepsilon, \bar{I}_t^\varepsilon)$, that is the solution to the Martingale problem corresponding to the following operator

$$A_\varepsilon f(x,i) := A_\varepsilon^{(i)} f(\cdot, i)(x) + \sum_{j=1}^J r_{ij}(x, x/\varepsilon) \left[ f(x,j) - f(x,i) \right] \tag{5.1}$$

with $A_\varepsilon^{(i)} : \mathcal{D}(A_\varepsilon^{(i)}) \subseteq C(\mathbb{R}^d) \to C(\mathbb{R}^d)$ be the generator of a strong $\mathbb{R}^d$-valued Markov process, with domain $\mathcal{D}(A_\varepsilon^{(i)})$.

## 5.1. Assumptions

Here we give the assumptions needed.

**Assumption 5.1.** The nonlinear generators $H_\varepsilon f = \varepsilon e^{-f/\varepsilon} A_\varepsilon e^{f/\varepsilon}$ admits an extended limit $H \subseteq ex - LIM H_\varepsilon$ with $H$ of the type

$$H := \left\{ (f, H_{f,\varphi}) \; : \; f \in C^2(\mathbb{R}^d), \; H_{f,\varphi} \in C(\mathbb{R}^d \times E') \text{ and } \varphi \in C^2(E') \right\}.$$

For all $\varphi$ there exist a map $\mathbf{H}_\varphi : \mathbb{R}^d \times \mathbb{R}^d \times E' \to \mathbb{R}$ such that for all $f \in D(H)$, $x \in \mathbb{R}^d$ and $z \in E'$, $H_{\varphi,f}(x,z) = \mathbf{H}_\varphi(x, \nabla f, z)$. Moreover, the image $\mathbf{H}_\varphi$ has the representation

$$\mathbf{H}_\varphi(x, p, z) = e^{-\varphi(z)} \left[ B_{x,p} + V_{x,p} + R_x \right] e^{\varphi}(z)$$

with $p = \nabla f(x)$ and $B_{x,p}, V_{x,p}, R_x$ such that

(i) For all $p \in \mathbb{R}^d$ there exists an eigenfunction $g_{x,p} \in C^2(\mathbb{R}^d \times J)$ with $g_{x,p}^i > 0$ and an eigenvalue $\mathcal{H}(x,p)$ such that $\left[ B_{x,p} + V_{x,p} + R_x \right] g_{x,p} = \mathcal{H}(x,p) g_{x,p}$.
(ii) $T_{x,p} = B_{x,p} + R_x$ verifies the maximum principle :
  if $(i_m, y_m) = \text{argmin } \varphi$ then $e^{-\varphi(i_m, y_m)} T_{x,p} e^{\varphi(i_m, y_m)} \geq 0$.
(iii) $p \mapsto V_{x,p}$ is coercive uniformly with respect to $x$.
(iv) $p \mapsto B_{x,p}$ and $p \mapsto V_{x,p}$ are convex uniformly on compact sets.

The above assumption implies the convergence of the nonlinear operators and the existence of the principal eigenvalue $\mathcal{H}$. Moreover, it will imply convexity and coercivity of $\mathcal{H}(x,p)$.

**Assumption 5.2.** The eigenvalue $\mathcal{H}$ is of the type $\mathcal{H}(x,p) = \sup_{\mu \in \mathcal{P}(E')} \left[ \Lambda(x,p,\mu) - I_{x,p}(\mu) \right]$ with

$$\Lambda(x,p,\mu) = \int_{E'} V_{x,p} \, d\mu, \quad \text{and} \quad I_{x,p}(\mu) = -\inf_{u>0} \int_{E'} \frac{(B_{x,p} + R_x)u}{u} \, d\mu,$$

and the following properties hold

(i) $I_{x,p}$ satisfies the assumption of Lemma 4.9,
(ii) $\Lambda(x,p,\mu)$ is continuous and $\|\Lambda(x,p,\mu)\|_\Theta < \infty$,
(iii) there exists a containment function $\Upsilon$ for $\Lambda$ in the sense of Definition 4.12,
(iv) for all $x$, there exists a unique measure $\mu_x^*$ such that $I_{x,0}(\mu_x^*) = 0$.

Assumption 5.2 implies the continuity of $\mathcal{H}$.

## 5.2. Large deviation for a switching Markov process

We are ready to state the general theorem.

**Theorem 5.3** (*Large Deviation for a Switching Markov Process*). *Let $(Y_t^\varepsilon, \bar{I}_t^\varepsilon)$ be the solution of the Martingale problem corresponding to the operator given in (5.1). If Assumptions 5.1 and 5.2 hold and suppose further that at time zero, the family of random variables $\{Y^\varepsilon(0)\}_{\varepsilon > 0}$ satisfies a large deviation principle in $\mathbb{R}^d$ with good rate function $\mathcal{I}_0 : \mathbb{R}^d \to [0, \infty]$. Then, the spatial component $\{Y_t^\varepsilon\}$ satisfies a large deviation principle in $C_{\mathbb{R}^d}[0, \infty)$.*

The proof of the above theorem follows the same lines of what is done in Section 4.3.

## 6. Conclusions and comparison with previous works

To conclude our work, in the following we summarise all the main novelties of our results.

1. We prove large deviations principle for the Markov process defined in Definition 2.1. The main steps of the proof are:
   (a) Convergence of the nonlinear generators $H_\varepsilon$ (3.2).
   (b) Proof of continuity of the Hamiltonian $\mathcal{H}$.
   (c) Comparison principle for $(1 - \lambda H)u = h$.
   (d) Proof of exponential tightness for $X^\varepsilon$.
   (e) Proof of the integral representation of the rate function.
2. We prove the Law of large numbers for the path characterising the limit process by calculating its speed. To do so, we also prove existence and uniqueness of the stationary measure of the operator (2.3).
3. We give a general result for a Switching Markov process.

The first result can be seen as an extension of part of the work [26]. We elaborate on the primary distinctions between our work and the previously cited one and how these distinctions contribute to the increased complexity of the proof of large deviations.

- We work on $\mathbb{R}^d$ and not on the torus $\mathbb{T}^d$. This transition from a compact to a non-compact set leads to the following significant consequences:
  (i) Firstly, in order to prove comparison principle, i.e step 1c, we need to construct four Hamiltonians in terms of a containment function that allows us to restrict part of our analysis on compact sets. Hence, we need to prove diagram 3. In [26], they only need two additional operators defined as multivalued limit of the Hamiltonian $H$.
  (ii) Secondly, in a compact setting step 1d is trivial. Indeed, the exponential tightness is implied by exponential compact containment condition that is always verified in a compact set.
  (iii) In the proof of the integral representation of the rate function, step 1e, some details are not needed in a compact setting as part of condition 8.9 in [14].

- We introduce a spatial component $x$ in the rates of the process that forces us to work with a Hamiltonian depending on both variables $x$ and $p$. For this reason, we need to work with a spatially inhomogeneous Hamilton–Jacobi equation. In particular:
  (i) Proving comparison principle one usually wants to bound the difference between sub-solution and super-solution $\sup_E u_1 - u_2$ by using a doubling variables procedure and typically ends up with an estimate of the following type
  $$\sup_E (u_1 - u_2) \leq \lambda \liminf_{\varepsilon \to 0} \left[ \mathcal{H}(x_\varepsilon, \alpha(x_\varepsilon - y_\varepsilon)) - \mathcal{H}(y_\varepsilon, \alpha(x_\varepsilon - y_\varepsilon)) \right] \qquad (6.1)$$
  $$+ \sup_E (h_1 - h_2).$$
  If the Hamiltonian does not depend on $x$, the final estimate is
  $$\sup_E (u_1 - u_2) \leq \lambda \liminf_{\varepsilon \to 0} \left[ \mathcal{H}(\alpha(x_\varepsilon - y_\varepsilon)) - \mathcal{H}(\alpha(x_\varepsilon - y_\varepsilon)) \right] + \sup_E (h_1 - h_2)$$
  $$= \sup_E (h_1 - h_2),$$
  that gives immediately comparison principle. This means that step 1c is partially immediate.
  (ii) To bound (6.1) in the non-spatially homogeneous case, we instead need to prove some regularity of the Hamiltonian. In particular, we need to prove continuity of $\mathcal{H}$ using some notions as $\Gamma$-convergence. If the Hamiltonian depends only on $p$, continuity, that is step 1b, follows immediately from convexity in $p$.

**Declaration of competing interest**

The authors declare that they have no known competing financial interests or personal relationships that could have appeared to influence the work reported in this paper.

**Acknowledgments**

This work was supported by The Netherlands Organisation for Scientific Research (NWO), grant number 613.009.148.

**Appendix. Uniqueness of the stationary measure**

We give here the proof of the existence and uniqueness of the stationary measure.

**Proposition Appendix A.1.** *Under the assumptions of Theorem 2.8, there exists a unique stationary measure of the operator*

$$L_{x,p}u(z,i) = \frac{1}{2}\Delta_z u(z,i) + (p - \nabla_x \psi^i(x,z)) \cdot \nabla_z u(z,i) + \sum_{j=1}^{J} r_{ij}(x,z)\left[u(z,j) - u(z,i)\right].$$

**Proof.** First of all, note that $L_{x,p}$ is of the type $T_{x,p} + R_x$ where $T_{x,p}$ is a diagonal matrix with diagonal elements $(T_{x,p})_{ii}$ of the type $\Delta + p \cdot \nabla + c$ and $(R_x)_{ii} = \sum_{j=1}^{J} r_{ij}$.

Let us consider, for some $\delta \in \mathbb{R}$, the operator $T_{x,p} + R_x + \delta$. For the latter operator, Conditions 1, 2 and 3 of [32] hold. Then, by [32, Theorem 1.1], there exists a unique function $\Psi \gg 0$ such that $(T_{x,p} + R_x + \delta)\Psi = \lambda\Psi$ for some $\lambda > 0$. It follows that $\lambda = \delta$ and $\Psi$ is the identity function. Hence, $\ker(T_{x,p} + R_x)$ is one-dimensional and it is spanned by the identity function, i.e., it consists of constants. Let $P_t$ be the semigroup associated to the generator $L_{x,p}$. By [11, Corollary V.4.6], $P_t$ is mean ergodic, that means that the Cesàro mean

$$C(r) = \frac{1}{r}\int_0^r P_s \, ds,$$

has a limit $P : C_b(E') \to C_b(E')$ for $r \to \infty$. Moreover, by [11, Lemma V.4.2], $Pf \in \ker(L_{x,p})$ for every $f \in C_b(E')$.

Let $T : c\mathbb{1} \in \ker(L_{x,p}) \mapsto c \in \mathbb{R}$. Then, $T \circ P : C_b(E') \to \mathbb{R}$ is a linear continuous function on $C_b(E')$. Then, by Riesz-Representation theorem, there exists a unique measure $\mu$ such that $(T \circ P)f = \langle f, \mu \rangle$. We show now that $\mu$ is the unique invariant measure for $P_t$. For all $f \in C_b(E')$ we have

$$\langle P_t f, \mu \rangle = (T \circ P)(P_t f) = (T \circ P)f = \langle f, \mu \rangle,$$

where in the second equality we used that $P = P_t P = PP_t = P^2$ (see [11, Lemma V.4.4]). Moreover, if $\mu^*$ is an invariant probability measure for $P_t$, let $Q$ be the projection $Q : f \in C_b(E') \mapsto \langle f, \mu^* \rangle \in \ker(L_{x,p})$. We show that $P = QP = Q$, obtaining then the uniqueness. On one hand, for $f \in C_b(E')$

$$QPf = \langle Pf, \mu^* \rangle = \langle \langle f, \mu \rangle, \mu^* \rangle = \langle f, \mu \rangle = Pf.$$

On the other hand,

$$QPf = \lim_{r \to \infty} \langle C(r)f, \mu^* \rangle = \lim_{r \to \infty} \langle f, C^*(r)\mu^* \rangle = \langle f, \mu^* \rangle = Qf. \quad \square$$